\documentclass[final]{siamltex}
\usepackage{amsfonts,amsmath,mathrsfs,bm}
\usepackage{amssymb}
\usepackage[notref,notcite]{showkeys}
\usepackage{hyperref}
\usepackage{graphicx}
\usepackage{wrapfig}
\usepackage{subfig}
\usepackage{float}
\usepackage{bbm}
\usepackage{algorithmic, algorithm}
\usepackage{gensymb,color,soul}
\usepackage{tikz}
\usepackage{tikz-3dplot}

\newcommand{\R}{\mathbb{R}}

\newcommand{\E}{\mathbb{E}}

\newtheorem{remark}[theorem]{Remark}

\DeclareMathOperator{\argmin}{\rm arg\, min}
\DeclareMathOperator{\argmax}{\rm arg\, max}

\title{Bayesian design of measurements for magnetorelaxometry imaging}

\author{
T. Helin\footnotemark[2]
\and N. Hyv\"onen\footnotemark[3]
\and J. Maaninen\footnotemark[3]
\and J.-P. Puska\footnotemark[3]
}

\begin{document}
\maketitle

\renewcommand{\thefootnote}{\fnsymbol{footnote}}
\footnotetext[2]{LUT University, School of Engineering Science, P.O.~Box 20, FI-53851 Lappeenranta, Finland. The work of TH was supported by the the Academy of Finland (decisions 326961, 345720, 348503, 353094).}
\footnotetext[3]{Aalto University, Department of Mathematics and Systems Analysis, P.O.~Box 11100, FI-00076 Aalto, Finland (nuutti.hyvonen@aalto.fi, jarno.maaninen@aalto.fi, juha-pekka.puska@aalto.fi). The work of NH, JM, and JP was supported by the Academy of Finland (decision 348503, 353081).}

\begin{abstract}
The aim of magnetorelaxometry imaging is to determine the distribution of magnetic nanoparticles inside a subject by measuring the relaxation of the superposition magnetic field generated by the nanoparticles after they have first been aligned using an external activation magnetic field that has subsequently been switched off.
    This work applies techniques of Bayesian optimal experimental design to (sequentially) selecting the positions for the activation coil in order to increase the value of data and enable more accurate reconstructions in a simplified measurement setup. Both Gaussian and total variation prior models are considered for the distribution of the nanoparticles. The former allows simultaneous offline computation of optimized designs for multiple consecutive activations, while the latter introduces adaptability into the algorithm by using previously measured data in choosing the position of the next activation. The total variation prior has a desirable edge-enhancing characteristic, but with the downside that the computationally attractive Gaussian form of the posterior density is lost. To overcome this challenge, the lagged diffusivity iteration is used to provide an approximate Gaussian posterior model and allow the use of the standard Bayesian A- and D-optimality criteria for the total variation prior as well. Two-dimensional numerical experiments are performed on a few sample targets, with the conclusion that the optimized activation positions lead, in general, to better reconstructions than symmetric reference setups when the target distribution or region of interest are nonsymmetric in shape.
\end{abstract}

\renewcommand{\thefootnote}{\arabic{footnote}}



\pagestyle{myheadings}
\thispagestyle{plain}
\markboth{T. HELIN, N. HYV\"ONEN, J. MAANINEN, AND J.-P. PUSKA}{DESIGN OF MEASUREMENTS FOR MRXI}

\section{Introduction}
\label{sec:introduction}

This work considers Bayesian {\em optimal experimental design} (OED) for choosing positions and orientations of activation coils in a simplified model for {\em magnetorelaxometry imaging} (MRXI). We implement both an offline algorithm for simultaneous optimization of multiple consecutive activations assuming a Gaussian prior for the imaged {\em magnetic nanoparticle} (MNP) distribution and an adaptive algorithm based on a {\em total variation} (TV) prior and the {\em lagged diffusivity iteration} as introduced for sequential X-ray imaging in~\cite{Helin22}.

\subsection{Magnetorelaxometry imaging}
In MRXI, the goal is to determine the distribution of magnetic nanoparticles inside a physical body from measurements on the relaxation of a superposition magnetic field generated by an alignment of the magnetic moments of the MNPs. The MNPs are aligned using an external activation magnetic field, which is then switched off, and the change in the magnetic field due to the relaxation of the MNPs is finally measured at sensors outside the examined body.

MNPs have a diameter of a few nanometers. They consist of ferro- or ferrimagnetic material and can be manipulated with an external magnetic field \cite{kotitz1997squid}. MNPs are employed in biomedicine \cite{pankhurst2003applications}: they have applications in different types of cancer treatment, e.g., in magnetic drug targeting \cite{torchilin00} and magnetic hyperthermia \cite{moroz2002tumor}. In these applications, it is important to have an estimate on the MNP distribution inside the subject because it directly affects how the treatment transpires.

As in \cite{Focke18}, we approximately model the dependence of the MRXI measurements of the relaxation magnetic field on the distribution of the MNPs via a linear forward model by working in the linear range of the Langevin function that is used for expressing the magnetization of the MNPs exposed to an external magnetic field. As further simplifications, we model the coil used for producing the external activation fields as a dipole and assume that each measurement sensor records the strength of the relaxation magnetic field (exactly) at a given point in a given direction. The  consecutive positions and orientations of the activation dipole are the design parameters we aim to optimize. 

After discretization, the forward model corresponds to a measurement matrix that depends nonlinearly on the design parameters, with a possibility to explicitly calculate the derivatives of this nonlinear dependence, which facilitates implementation of differentiation-based optimization algorithms. Due to a certain symmetry in the measurement model, the system matrix depends essentially in the same way on the measurement positions and directions as on the specifications of the activation dipoles, which means that the introduced techniques could as well be used for designing the configuration of the measurement sensors.

\subsection{Bayesian experimental design}

A Bayesian optimal design $p^*$ is defined as a maximizer over the design space $\mathcal{P}$ for the expectation of the utility function $\E_{u,y}[U( p ; u,y)]$ with respect to the data $y \in \mathcal{Y}$ and the model parameters $u \in \mathcal{U}$ \cite{chaloner1995bayesian}. That is,
\begin{equation}
	\label{eq:OED_task}
	p^* 
	 =  \underset{p\in \mathcal{P}}{\argmax} \int_{\mathcal{Y}} \int_\mathcal{U} U(p; u, y) \pi(u \, | \, p, y) \pi(y \, | \, p) \, {\rm d} u \, {\rm d}y,
\end{equation}
where $\pi(u \, | \, p, y)$ and $\pi(y \, | \, p)$ are the posterior distribution of the parameter $u$ and the marginal distribution of the data $y$, respectively, under the design $p$. Two of the most common choices for $U$ are a \emph{negative quadratic loss function} that measures the squared distance from $u$ to a specific point estimator such as the posterior mean and the \emph{expected information gain} for which $U$ is related to the Kullback--Leibler distance between the posterior and prior distributions. 

In our setting of MRXI, the design parameter vector $p$ encodes the consecutive positions and orientations of the dipole-like activation coil, $u$ corresponds to the discretized MNP distribution, and $y$ is the vector of measurements on the magnetic field at the sensors. As the measurement model relating $u$ and $y$ is linear, assuming a Gaussian prior and an additive Gaussian noise model considerably simplifies the optimization targets corresponding to the  aforementioned two utility functions: employing a quadratic loss function or maximizing the information gain lead to minimizing a weighted trace or the determinant of the posterior covariance, respectively. These correspond to so-called A- and D-optimal designs \cite{alexanderian2016bayesian,chaloner1995bayesian}. What is more, the posterior covariance is independent of the measurements, meaning that the experimental design can be performed offline,~i.e.~prior to taking any measurements, and simultaneously for several consecutive activations.

If a smoothened TV prior is utilized for the MNP distribution, it is not possible to get rid of the double integral over the potentially high-dimensional spaces in \eqref{eq:OED_task} without further simplifications. However, if one proceeds sequentially, choosing the specifications for the next activation only after computing a {\em maximum a posteriori} (MAP) estimate for the MNP distribution via the lagged diffusivity iteration \cite{Vogel96} based on the measurements from the previous activations, it is possible to interpret the MAP estimate as the mean of a Gaussian distribution whose covariance matrix is available as a side product of the iteration \cite{Bardsley18,Calvetti08}. Basing the selection of the specifications for the next dipole activation on this covariance structure, one can devise a sequential Bayesian OED algorithm that adapts to the already collected data and has potential to produce edge-promoting experimental designs; see \cite{Helin22} for an application of this idea to sequential X-ray imaging. Let us also note that the sequential optimization approach can be extended to a non-parametric setting with convex priors (such as the smoothened TV) for which the posterior has good approximation properties by Gaussian distributions in the neighborhood of non-parametric MAP estimators~\cite{helin2015maximum}.

\subsection{Our contribution}
The main contribution of this work is the application of Bayesian OED to a simplified two-dimensional model of MRXI, where the activation coils and measurement sensors are modeled as point-like objects with orientations. OED has previously been studied in the framework of MRXI in \cite{Maaninen23,Schier21,Schier20}, of which only the master's thesis \cite{Maaninen23} considers a Bayesian setting. 

The presented numerical experiments tackle OED both with a Gaussian prior and simultaneous offline optimization of several activations and with a smoothened TV prior and sequential adaptive optimization of the activations. In both cases, the reconstruction accuracy for the optimal designs is compared to symmetric reference setups, demonstrating that the employment of Bayesian OED improves the performance of MRXI. However, as the the objective functions considered when deducing the `optimal' designs suffer from multiple local optima, and a greedy approach is used when optimizing the activations sequentially, there is no certainty that globally optimal designs are actually found in all experiments.

For the simultaneous design of activations, the optimization is performed by Newton's method or  gradient descent due to the high-dimensional design space, which makes the approach prone to finding local optima. In the sequential algorithm, where only the position and orientation of the next activation dipole need to be optimized at a time, an exhaustive search is also considered in order to evaluate the optimality of the designs produced by the differentiation-based methods. The numerical experiments focus mainly on A-optimality, but extending all presented considerations to the case of D-optimality would be conceptually straightforward.

The considered sequential approach to Bayesian OED has previously been investigated for choosing optimal projection geometries in X-ray imaging with a Gaussian prior in \cite{Burger21} and with a TV prior in \cite{Helin22}. However, these papers do not tackle simultaneous optimization of many projection geometries with a Gaussian prior due to computational restrictions. Moreover, compared to X-ray tomography with a limited projection aperture, the measurement matrix for a single activation in MRXI is more ill-conditioned and provides information about the entire imaged object. Hence, the functionality of the sequential OED algorithm with the smoothened TV prior for MRXI is not obvious based solely on the material in \cite{Helin22}. 

Our method relies on a well-defined discretization of an underlying non-parametric Bayesian OED problem. Non-parametric approach for OED in Bayesian inverse problems has been formalized by Alexanderian (see \cite{alexanderian2016bayesian} and the references therein), and it provides the theoretical underpinnings for many applications emerging,~e.g.,~in inverse problems related to PDEs (see,~e.g.,~\cite{alexanderian2016bayesian, alexanderian2014optimal,  beck2018fast, crestel2017optimal, long2015fast, wu2020fast}) or nonlinear systems \cite{huan2013simulation}. Stability properties of the expected utility concerning model approximations in Bayesian optimal experimental design were recently investigated in \cite{duong2022stability}. This finding suggests that, in our framework, each optimization step in the sequential D-optimal approach remains robust when subjected to discretization and linearization of the Langevin function. However, rigorous study of such robustness is part of future work. For more general review on Bayesian OED, see \cite{alexanderian2021optimal_review, chaloner1995bayesian, rainforth2023modern, ryan2016review}.

This text is organized as follows. The (simplified and discretized) measurement model of MRXI is described in Section~\ref{sec:mm}. Section~\ref{sec:priors} introduces the Bayesian framework for inverse problems and recalls the probabilistic interpretation of the lagged diffusivity iteration from~\cite{Bardsley18,Helin22}. In Section~\ref{sec:oed}, the optimality targets for A- and D-optimal designs are introduced, and their minimization is considered. Section~\ref{sec:opti} presents the complete optimization/reconstruction algorithm for the sequential design process, and Section \ref{sec:numerics} is dedicated to numerical experiments. Finally, Section~\ref{sec:conclusion} lists the concluding remarks.

\section{Measurement model and its discretization}
\label{sec:mm}

MRXI consists of two main phases that are excitation and relaxation. In the excitation phase, an external magnetic field is generated by using electromagnetic coils, called {\em activations} in what follows. The activations realign the magnetic moments of the MNPs so that the superposition field that they generate can be measured with magnetic field sensors outside the imaged object. In the relaxation phase, the activations are switched off, allowing the magnetic moments of the MNPs to reorient through N{\'e}el and Brownian relaxations until equilibrium is reached~\cite{kotitz1995time}. N{\'e}el relaxation describes the rotation of the magnetic moment within the core of an MNP and Brownian relaxation the rotation of the entire MNP. Note that magnetic field sensors can only measure a change in a magnetic field, meaning that the MRXI measurements actually correspond to the change in the field over some time period during the relaxation phase.

Ideally, a single MRXI measurement equals the strength of the superposition magnetic field created by the MNPs at a given location in a given direction. To achieve this, the external activation magnetic field should vanish immediately after it is switched off and the employed measurement sensor should be able to measure the change in the magnetic field precisely at a single point, both of which are impossible conditions to satisfy in practice. Be that as it may, in the following analysis we assume such idealized measurements and refer to \cite{Jaufenthaler20,liebl2014quantitative} and references therein for information on the practical limitations of MRXI instrumentation.

The imaged object is represented by a domain $\Omega \subset \mathbb{R}^d$, with $d=3$, and the density of the MNPs in $\Omega$ is modeled by $c \in L^2_{+}(\Omega)$, where 
$$
L^2_+(\Omega) = \big\{ v \in L^2(\Omega) \ | \ {\rm ess} \inf v \geq  0 \big\} 
$$
accounts for the presumed nonnegativity of the MNP density. The aim of MRXI measurements is to reconstruct $c$ based on measurements on (the relaxation of) superposition magnetic fields generated by the MNPs under a series of activations.

\subsection{Modeling the activations}
Consider an external activation coil at $a \in \mathbb{R}^d \setminus \Omega$, with its shape described by a smooth enough (closed) path $\Gamma_a \subset \R^d$. Assume furthermore that
\begin{equation*}
    \gamma_a: [0, L] \to \R^d
\end{equation*}
gives an arclength parametrization for $\Gamma_a$. According to the Biot--Savart law, the magnetic field induced when a constant net current $J$ runs through the coil is given as
\begin{equation}
\label{eq:biot_savart_coil}
B_a(x) =  \frac{\mu_0}{4\pi} J  \int_0^L \frac{\gamma_a'(\tau) \times (x - \gamma_a(\tau))}{|x - \gamma_a(\tau)|^3} {\rm d} \tau, \qquad x \in \R^d \setminus \Gamma_a,
\end{equation}
where $| \cdot |$ denotes the Euclidean norm.
A common approach to modeling electromagnetic coils is to approximate the coil path with linear segments. This approximation is explained in \cite{Hanson02} and used for an example on MRXI in \cite{liebl2014quantitative}. The issue with such a numerical implementation is the computational complexity of the MRXI forward model, especially when the coils are approximated to a high precision. Moreover, the detailed modeling of the coils is case-dependent, and it is not obvious that such details have a considerable effect on the general conclusions on the applicability of Bayesian OED to MRXI. For these reasons, we model the activations as electromagnetic dipoles; cf.~\cite{Focke18}.

Our use of dipole activations can be motivated, e.g., by assuming that $\Gamma_a$ consists of a single circular loop of radius $\rho$ in the $(d-1)$-dimensional plane that contains the center point of the loop $a$ and has the unit vector $\nu$ as its normal. If $\rho$ is small compared to $|x-a|$, it is well known that $B_a(x)$ can be approximated by the field of a magnetic dipole at $a$. More precisely (see,~e.g.,~\cite{Seleznyova16}),
\begin{equation}
\label{eq:asymp_dipole}
    B_a(x) = \pm \frac{\mu_0}{4\pi} \bigg( \frac{3(x-a)(x-a)^\top}{|x - a|^5} - \frac{I}{|x-a|^3} \bigg) \alpha + O\Big( \frac{\rho^4}{|x - a|^5} \Big),
\end{equation}
where $I \in \R^{d\times d}$ is the identity matrix, $\alpha = \pi \rho^2 J \nu$ has an interpretation as a dipole moment, and the sign depends on the direction of the current in the loop and the chosen orientation for $\nu$. In the following, we drop the latter term in \eqref{eq:asymp_dipole} and also simplify/abuse the notation by writing
\begin{equation*}
    B_a(x) = \bigg( \frac{3(x-a)(x-a)^\top}{|x - a|^5} - \frac{I}{|x-a|^3} \bigg) \alpha,\qquad x \in \R^d \setminus \{ a \},
\end{equation*}
where $\alpha \in \R^d$ is a dipole moment that includes all physical constants and will together with $a$ act as a  design parameter in our numerical experiments.

\subsection{Modeling the measurements}
Let $M: \Omega \to \R^d$ describe the magnetization generated by the MNPs inside the imaged object due to the external magnetic field. The superposition field generated by the magnetization is expressed as
\begin{equation*}
B_M(w) = \frac{\mu_0}{4\pi} \int_\Omega \bigg( \frac{3(w-x)(w-x)^\top}{|w - x|^5} - \frac{I}{|w-x|^3} \bigg) M(x) \, {\rm d} x, \qquad w \notin \overline{\Omega},
\end{equation*}
that is, as a field induced by a density of magnetic dipoles over $\Omega$ described by the magnetization. According to the basic theory on idealized paramagnetic materials \cite{Reitz09}, the magnetization $M$ generated during the excitation phase can be modeled as~(cf.~\cite{Focke18})
\begin{equation}
\label{eq:mnp_magnetization}
M(x) \propto \mathcal{L}\big( q B_a(x)\big)c(x),\qquad x \in \Omega,
\end{equation}
where $q>0$ is a MNP-dependent physical constant, $B_a$ is the activation magnetic field and $\mathcal{L}: \R \to \R$ is the Langevin function
\begin{equation}\label{eq:langevin}
\mathcal{L}(\tau) = \mathrm{coth}(\tau) - \frac{1}{\tau} = \frac{1}{3} \tau + O\big(|\tau|^3\big)
\end{equation}
applied in \eqref{eq:mnp_magnetization} componentwise.

Since in MRXI the argument of $\mathcal{L}$ in \eqref{eq:mnp_magnetization} can be assumed to be in the linear range of the Langevin function~\cite{Focke18}, we drop the second term on the right-hand side of \eqref{eq:langevin}. Assuming further that a measurement sensor at $s \in \mathbb{R}^d \backslash \Omega$  measures the (relaxation-induced change in the superposition) magnetic field in the direction of a vector $\sigma \in \R^d$ that describes the sensor's orientation, the corresponding measurement is modeled as
\begin{equation}
\label{eq:nondisc_measurement}
y_{s,\sigma} = \sigma^\top  B_M(s) =   \sigma^\top \int_\Omega \bigg( \frac{3(s-x)(s-x)^\top}{|s - x|^5} - \frac{I}{|s-x|^3} \bigg) B_a(x) \, c(x) \, {\rm d} x.
\end{equation}
Observe that in \eqref{eq:nondisc_measurement} all constants have been included in the measurement direction vector $\sigma$.

\subsection{Complete model and its discretization}

Combining our activation and measurement models and recalling the assumption to be able to measure the strength of the actual superposition magnetic field via measuring its relaxation, the complete measurement model for a single measurement due to a single activation reads
\begin{equation}
\label{eq:complete_nondisc_measurement}
y_{s,\sigma, a, \alpha} =  \big(\mathcal{K}_{\sigma, a,\alpha}c \big)(s) := \int_\Omega \kappa_{\sigma,a,\alpha}(s,x) c(x) \, {\rm d} x,
\end{equation}
where the integral kernel is
\begin{equation*}
    \kappa_{\sigma,a,\alpha}(w,x) = \sigma^\top \bigg(\frac{3 (w-x)(w-x)^\top}{|w - x|^5} - \frac{I}{|w-x|^3} \bigg) \bigg( \frac{3(x-a)(x-a)^\top}{|x - a|^5} - \frac{I}{|x-a|^3} \bigg) \alpha.
\end{equation*}
The integral operator $\mathcal{K}_{\sigma,s,a}: L^2(\Omega) \to L^2(D)$, implicitly defined by \eqref{eq:complete_nondisc_measurement}, is compact due to the boundedness (or smoothness) of its kernel, assuming the domain for the possible positions of the sensor locations $D$ is bounded and satisfies $\overline{D} \cap \overline{\Omega} = \emptyset$ \cite{Ringrose71}. This gives an explanation for the ill-posedness of the considered MRXI problem~\cite{Focke18}. Observe also that \eqref{eq:complete_nondisc_measurement} represents a linear dependence between the MNP concentration $c \in L^2_+(\Omega)$ and the measurement $y_{s,\sigma, a, \alpha} \in \R$, parametrized by the two vector pairs $(a, \alpha)$ and $(s, \sigma)$ that define the activation dipole and the measurement sensor, respectively. Moreover, it is easy to check via transposition that $\kappa_{\sigma,a,\alpha}(s, \, \cdot \,) = \kappa_{\alpha,s,\sigma}(a,\, \cdot \,)$, which means that changing the roles of the pairs $(a, \alpha)$ and $(s, \sigma)$ does not change the measurement,~i.e.,~$y_{s,\sigma,a,\alpha} = y_{a,\alpha,s,\sigma}$. In particular, our approach for designing the activation pattern introduced in what follows could as well be used for designing the positions of the measurement sensors.

The relation \eqref{eq:complete_nondisc_measurement} can be discretized as
\begin{equation}
\label{eq:disc_meas}
    y_{s,\sigma, a, \alpha} \approx \sum_{j=1}^{N_c} \big(k_{s,\sigma,a,\alpha}\big)_j \, c_j = k_{s,\sigma,a,\alpha}^\top c
\end{equation}
where $c_1, \dots, c_{N_c}$ are the degrees of freedom in the employed parametrization for the MNP concentration and we have abused the notation by denoting now with $c$ a vector in $\R^{N_c}$. The coefficient vector $k_{s,\sigma,a,\alpha}$ is defined by the function $\kappa_{\sigma,a,\alpha}(s, \, \cdot \,)$ and the quadrature rule utilized for numerically evaluating the right-hand side of  \eqref{eq:complete_nondisc_measurement}. In our numerical experiments, the two-dimensional\footnote{The interpretation of MRXI measurements in two spatial dimensions is considered in Section~\ref{sec:numerics}.} domain $\Omega$ is divided into a homogeneous grid of $N_c$ pixels with centers $x_1, \dots, x_{N_c}$, the degrees of freedom for the MNP concentration are $c_j = c(x_j)$, $j=1, \dots, N_c$, and
\begin{equation*}
    \big(k_{s,\sigma,a,\alpha}\big)_j = \omega \,\kappa_{\sigma,a,\alpha}(s, x_j), \qquad j = 1, \dots, N_c,
\end{equation*}
where the weight $\omega$ is the area of a single pixel.

Let us assume that there are $N_s$ measurements for each activation; the corresponding parameter pairs $(s_j, \sigma_j)$, $j=1, \dots, N_s$, are not considered as design parameters, but they are predefined and the same for each activation. According to \eqref{eq:disc_meas}, a full set of measurements for a single activation can thus be modeled as
\begin{equation}
    y_{p} = K(p) c \in \R^{N_s},
\end{equation}
where $p = (a, \alpha)$ denotes the design parameter pair,~i.e.,~the position and moment for the activation dipole, and
\begin{equation*}
    K(p) = \begin{bmatrix}
    k_{s_1, \sigma_1, a, \alpha}^\top \\
    \vdots \\
    k_{s_{N_s}, \sigma_{N_s}, a, \alpha}^\top
    \end{bmatrix} \in \R^{N_s \times N_c}.
\end{equation*}
Analogously, after $k$ sets of measurements corresponding to the activations defined by a sequence of design parameters $p_j = (a_j, \alpha_j)$, $j=1, \dots,k$, the discretized noiseless measurement model reads
\begin{equation}
\label{eq:ax_final}
\mathbf{y}_k = \begin{bmatrix}y_1 \\ \vdots \\ y_k \end{bmatrix} = \mathbf{K}(\mathbf{p}_k)  c, \qquad \text{with} \quad \mathbf{K}(\mathbf{p}_k) = \begin{bmatrix}K(p_1) \\ \vdots \\ K(p_{k}) \end{bmatrix} \in \R^{kN_s \times N_c},
\end{equation}
where $\mathbf{p}_k$ denotes a concatenation of all $k$ design variables. This is the model that is used for (sequential) Bayesian OED in the following. In particular, we do not insist on the nonnegativity of the discretized MNP concentration in what follows.

\section{Prior models and lagged diffusivity iteration}
\label{sec:priors}

In this section, a Bayesian model for the measurements deterministically described by \eqref{eq:ax_final} is introduced and two different prior models for the discretized MNP concentration are considered. The first prior model is Gaussian, which enables optimization of the measurements in an offline mode and, in particular, simultaneous optimization of the specifications of several sequential activations. The second one is a smoothened TV prior, which makes sequential optimization of the activations truly adaptive, that is, the previously collected measurement data affect the subsequent activation designs.

Consider the probabilistic measurement model for a single activation,~i.e.,~for a single patch of $N_s$ rows in \eqref{eq:ax_final},
\begin{equation}
  \label{eq:meas_model}
Y_j = K(p_j) C + E_j, \qquad j=1, \dots, k,
\end{equation}
where $C$ is the discretized and MNP concentration now modeled as a random vector, and the noise $E_j$ is assumed to follow a zero-mean Gaussian distribution $\mathcal{N}(0, \Gamma_{{\rm noise}, j})$, where is $\Gamma_{{\rm noise}, j} \in \R^{N_s \times N_s}$ is symmetric and positive definite.  The maximum number of activations is denoted by $N_a$, meaning that $1 \leq k \leq N_a$ in all following considerations. The noise processes $E_1, \dots, E_{N_a}$ are assumed to be mutually independent.

The prior for the MNP density $C$ is assumed to be independent of the noise processes and to follow a probability density of the form
\begin{equation}
\label{eq:sigma_prior}
\pi(c) \propto \exp \! \big( -\gamma  \Phi(c) \big),
\end{equation}
where the role of $\gamma > 0$ is separately specified for the two considered priors. According to the Bayes' formula and assuming the measurement model \eqref{eq:meas_model}, the posterior density for $c$ after $k$ activations is
\begin{align}
  \label{eq:post_u}
\pi(c \, | \, \mathbf{y}_k) \, &\propto \, \pi( \mathbf{y}_k \, | \, c) \, \pi(c)  \nonumber \\
&\propto \, \exp \Big(-\frac{1}{2}  (\mathbf{y}_k - \mathbf{K}(\mathbf{p}_k) c )^{\rm T} (\bm{\Gamma}_{{\rm noise}, k})^{-1}(\mathbf{y}_k - \mathbf{K}(\mathbf{p}_k) c ) - \gamma \Phi(c) \Big),
\end{align}
where $\bm{\Gamma}_{{\rm noise}, k} := {\rm diag}(\Gamma_{{\rm noise}, 1}, \dots, \Gamma_{{\rm noise}, k})\in \R^{kN_s \times kN_s}$ is a block diagonal matrix defined by the noise covariance matrices for the previous measurements.

\subsection{Gaussian prior}
Let us first assume that {\em a priori} $C \sim \mathcal{N}(\widehat{c}_0, \Gamma_0)$, where $\widehat{c}_0 \in \R^{N_c}$ is the prior mean and the symmetric positive-definite $\Gamma_0 \in \R^{N_c \times N_c}$ is the prior covariance. Hence, $\gamma = 1/2$ and
\begin{equation*}
\Phi(c) = (c - \widehat{c}_0)^\top \Gamma_0^{-1}(c - \widehat{c}_0)
\end{equation*}
in \eqref{eq:sigma_prior}. For a Gaussian prior, it is well known that the posterior $\pi(c \, | \, \mathbf{y}_k)$ is also a Gaussian (see,~e.g.,~\cite{Kaipio06}), with the covariance
  \begin{align}
  \label{eq:Gaussian_posterior}
    \Gamma_k &=  \Gamma_0 - \Gamma_0 \mathbf{K}(\mathbf{p}_k)^\top\big(\mathbf{K}(\mathbf{p}_k) \Gamma_0 \mathbf{K}(\mathbf{p}_k)^\top + \bm{\Gamma}_{{\rm noise},k} \big)^{-1} \mathbf{K}(\mathbf{p}_k) \Gamma_0 \\
    &  =\Gamma_{k-1} - \Gamma_{k-1} K(p_k)^\top\big(K(p_k) \Gamma_{k-1} K(p_k)^\top + \Gamma_{{\rm noise},k} \big)^{-1} K(p_k) \Gamma_{k-1}
    \label{eq:Gaussian_posterior_one_step}
   \end{align}
and the mean 
\begin{align*}
  \widehat{c}_k &= \widehat{c}_0 + \Gamma_0 \mathbf{K}(\mathbf{p}_k)^\top\big(\mathbf{K}(\mathbf{p}_k) \Gamma_0 \mathbf{K}(\mathbf{p}_k)^\top + \bm{\Gamma}_{{\rm noise},k} \big)^{-1}\big(\mathbf{y}_k - \mathbf{K}(\mathbf{p}_k)^\top \widehat{c}_0\big) \\
  &= \widehat{c}_{k-1} +  \Gamma_{k-1} K(p_k)^\top\big(K(p_k) \Gamma_{k-1} K(p_k)^\top + \Gamma_{{\rm noise},k} \big)^{-1} (y_k - K(p_k) \widehat{c}_{k-1}),
  \end{align*}
    where the latter recursive formulas follow by treating the posterior after the previous measurement as the prior for the newest one.

It is important to notice that the posterior covariance after $k$ measurements depends on the previous experimental designs via $\mathbf{K}(\mathbf{p}_k)$ (cf.~\eqref{eq:ax_final}), but it does {\em not} depend on the corresponding measurement data. As the optimization targets considered in Section~\ref{sec:oed} are functions of the posterior covariance only, the optimal designs do not depend on measured data either. Hence, the optimization of the measurement setup can be performed offline, with no other reason to resort to greedy sequential optimization of the activations than the computational cost related to considering a high-dimensional decision variable. As a consequence, we aim to simultaneously optimize all activation designs $\mathbf{p}_{N_a}$  when considering a Gaussian prior for the MNP concentration.

\subsection{Total variation prior and lagged diffusivity}
\label{sec:TV}

For smoothened TV, the scaled negative log-prior $\Phi$ is defined as 
\begin{equation}
\label{eq:aRRa}
\Phi(c) = \int_{\Omega} \varphi \big(|\nabla c | \big) \, {\rm d} x, \qquad \text{with} \quad \varphi(t) = \sqrt{t^2 + T^2} \approx | t |,
\end{equation}
accompanied by the information that $c$ vanishes at the pixels next to the boundary of~$\Omega$. The small parameter $T>0$ ensures the differentiability of $\varphi$, which is required by the lagged diffusivity iteration, and in this case $\gamma>0$ in \eqref{eq:sigma_prior} controls the strength of the prior. With this choice, the posterior \eqref{eq:post_u} is obviously not Gaussian, which typically makes the evaluation of  an optimality target of the form \eqref{eq:OED_task} computationally demanding. To circumvent this problem, we adopt the approach in~\cite{Helin22} and iteratively approximate $\Phi(c)$ by quadratic terms in the spirit of the lagged diffusivity iteration~\cite{Vogel96}; see also~\cite{Arridge13,Harhanen15}. With a suitable stopping condition, this technique automatically produces a Gaussian approximation for the posterior of $c$ after $k$ activations; the covariance matrix for this approximate posterior can then be deployed when determining the experimental design for the next activation. In particular, this approach leads to a sequential OED algorithm that adapts to the measurement data in hand.

To set the stage for introducing a Bayesian version of the lagged diffusivity iteration, let us identify $c$ with its interpolant in a piecewise linear finite element basis $\{ \phi_j\}_{j=1}^{N_c} \subset H^1_0(\Omega)$. A straightforward differentiation reveals that
$$
\nabla_{c} \Phi(c) = \Theta(c) c, \qquad c \in \R^{N_c},
$$
where
\begin{align}
\label{eq:H}
\Theta_{i,j}(c)
&:= \int_{\Omega} \frac{1}{\sqrt{|\nabla_x c(x)|^2 + T^2}} \, \nabla \phi_i(x) \cdot \nabla \phi_j(x) \, {\rm d} x,
\qquad i,j=1,\dots, N_c,
\end{align}
for any $c \in \R^{N_c}$ interpreted as an element of $H^1_0(\Omega)$ via the introduced finite element basis. In particular, $\Theta$ is positive definite and thus invertible since it corresponds to a finite element discretization of an elliptic partial differential operator accompanied by a homogeneous Dirichlet boundary condition; see~\cite{Helin22} for more details.

Let us then assume that we have been able to deduce (an approximation of) the MAP estimate $\widehat{c}_{k-1}$ for the posterior \eqref{eq:post_u} after $k-1$ activations, with $\Phi$ defined by \eqref{eq:aRRa}. If $k=1$, our initial guess for the MNP concentration plays the role of $\widehat{c}_{k-1}$. The aim is to use the lagged diffusivity iteration to compute the MAP estimate $\widehat{c}_{k}$ for the posterior \eqref{eq:post_u} after $k$ activations in such a way that we simultaneously form a Gaussian approximation for \eqref{eq:post_u}.
The lagged diffusivity iteration is started by setting $\widehat{c}_{k}^{(0)} = \widehat{c}_{k-1}$. The subsequent iterates are defined recursively via
\begin{equation}
  \label{eq:step2}
 \widehat{c}_{k}^{(j)} =  \Gamma_{k}^{(j-1)} \mathbf{K}(\mathbf{p}_k)^\top \big(\mathbf{K}(\mathbf{p}_k) \Gamma_{k}^{(j-1)} \mathbf{K}(\mathbf{p}_k)^\top + \gamma \bm{\Gamma}_{{\rm noise}, k}  \big)^{-1}\mathbf{y}_k,
\end{equation}
where $\Gamma_{k}^{(j-1)} = \Theta(\widehat{c}_{k}^{(j-1)})^{-1}$ can be interpreted as the covariance matrix for a zero-mean Gaussian prior multiplied by $\gamma$. The iterate itself $\widehat{c}_{k}^{(j)}$ has an interpretation as the corresponding posterior mean after performing measurements corresponding to $\mathbf{K}(\mathbf{p}_k)$. Although not needed explicitly in the lagged diffusivity iteration itself, the posterior mean can be accompanied by a Gaussian density with the covariance
\begin{equation}
\label{eq:posterior_cov}
    \Gamma_{k}^{(j)} = \gamma^{-1} \big( \Gamma_{k}^{(j-1)} - \Gamma_{k}^{(j-1)} \mathbf{K}(\mathbf{p}_k)^\top \big(\mathbf{K}(\mathbf{p}_k) \Gamma_{k}^{(j-1)} \mathbf{K}(\mathbf{p}_k)^\top + \gamma \bm{\Gamma}_{{\rm noise}, k}  \big)^{-1} \mathbf{K}(\mathbf{p}_k) \Gamma_{k}^{(j-1)} \big);
\end{equation}
see once again \cite{Helin22} for more details. 

This iterative process is continued until a suitable stopping criterion is satisfied, say, at $j = J$; the criterion employed in our numerical experiments is considered in Section~\ref{sec:opti}. One then defines $\widehat{c}_k = \widehat{c}_{k}^{(J)}$ to be the reconstruction after $k$ projection images,~i.e.,~an approximation of the MAP estimate for \eqref{eq:post_u}. Moreover, the covariance matrix for the associated approximate Gaussian density $\Gamma_k = \Gamma_{k}^{(J)}$ is used for choosing the design parameter vector $p_{k+1}$ defining the next activation. In particular, the covariance matrix \eqref{eq:posterior_cov} only needs to be evaluated once at $j=J$ when the iteration is terminated. For more details on the lagged diffusivity iteration see \cite{chan1999convergence,dobson1997convergence,Vogel96}.

\section{Computing optimal designs}
\label{sec:oed}
In this section, we recall the concepts of A- and D-optimal designs and consider numerically solving the associated optimization problems. The presentation is intentionally compact; we refer to~\cite{alexanderian2016bayesian,chaloner1995bayesian} and \cite{Nocedal06}, respectively, for more information on the optimality conditions of Bayesian OED and the tools of nonlinear optimization.

\subsection{A- and D-optimality}
\label{ssec:A_optimality}

The A- and D-optimality criteria aim to minimize the (weighted) trace and the determinant of the (Gaussian) posterior, respectively, with respect to the design parameters. The design parameter vector is denoted by $\xi$ and the considered posterior covariance by $\Sigma(\xi)$. These entities can have two meanings corresponding to two different settings:
\begin{enumerate}
\item The design $\xi$ parametrizes $\mathbf{p}_{N_a}$,~i.e.,~it defines the specifications of all activations, and $\Sigma(\xi) = \Gamma_{N_a}(\mathbf{p}_{N_a}(\xi))$ is the final posterior for a Gaussian prior defined by \eqref{eq:Gaussian_posterior} with $k = N_a$.
\item The measurements corresponding to the first $k-1$ activations have already been performed, $\xi$ parametrizes the $k$th activation, and $\Sigma(\xi)$ is defined in accordance with \eqref{eq:Gaussian_posterior_one_step} as
\begin{equation}
\label{eq:post_seq}
 \Sigma(\xi) = \Gamma_{k-1} - \Gamma_{k-1} K(p)^\top\big(K(p) \Gamma_{k-1} K(p)^\top + \Gamma_{{\rm noise},k} \big)^{-1} K(p) \Gamma_{k-1},
\end{equation}
where $p = p(\xi)$ and $\Gamma_{k-1}$ is the covariance of a Gaussian distribution that the MNP concentration is assumed to follow after $k-1$ measurements.
\end{enumerate}
The main motivation for considering the latter case is the method for sequentially building approximate Gaussian posteriors under a smoothened TV prior reviewed in Section~\ref{sec:TV}, but it can in principle also be used for deducing greedy sequential designs under a Gaussian prior (cf.,~e.g.,~\cite{Burger21}).

In our setting, an A-optimal design is defined as
\begin{equation}
\label{eq:Aoptimal}
\xi_{\rm A} =  \underset{\xi}{\argmin} \, \, \Psi_{\rm A}(\xi), \qquad \text{with} \quad \Psi_{\rm A}(\xi) = {\rm tr}  \big(A \Sigma(\xi) A^\top\big).
\end{equation}
It minimizes the expected squared distance of the unknown from the posterior mean in the seminorm defined by the positive semidefinite matrix $A^\top \!A$ for a given $A \in \R^{N_c \times N_c}$;~see, e.g., \cite{alexanderian2016bayesian,chaloner1995bayesian} for more details. In our numerical tests, $A$ is usually the identity matrix~$I$, indicating that all degrees of freedom in the MNP concentration are considered equally important. However, certain diagonal elements of the identity matrix can also be set to zero in order to only account for the reconstruction error in the other degrees of freedom. Such a diagonal matrix is denoted by $I_{\rm ROI}$, where "ROI" indicates that the reconstruction error over some {\em region of interest} (ROI) is minimized. In our numerical tests, where $\Omega$ is divided into $N_c$ homogeneous pixels with constant concentration values, an estimate of the expected $L^2$-error over the ROI for a design $\xi$ can be computed as
\begin{equation}
\label{eq:L2_error}
\tilde{\Psi}_{\rm A}(\xi) = \sqrt{\frac{|\Omega|}{N_c} \Psi_{\rm A}(\xi)},
\end{equation}
where $A = I_{\rm ROI}$ and $|\Omega|$ denotes the area of $\Omega$.

A D-optimal design maximizes the information gain when the prior is replaced by the posterior~\cite{chaloner1995bayesian}, which in our setting can be expressed as
\begin{equation}
\label{eq:Doptimal}
\xi_{\rm D} =  \underset{\xi}{\argmin} \, \Psi_{\rm D}(\xi), \qquad \text{with} \quad \Psi_{\rm D}(\xi) = \log ( \det \Sigma(\xi)).
\end{equation}
The minimization target $\Psi_{\rm D}$ equals the negative of the information gain up to an additive constant and scaling, and the inclusion of a logarithm in $\Psi_{\rm D}$ follows from information theory, but it also makes the evaluation of $\Psi_{\rm D}$ more stable (cf.~\cite{Burger21}). It would also be possible to only consider information gain over some ROI \cite{Burger21}, but such a case is not considered in the numerical examples of this work. 

In the rest of this section, the minimization target is denoted generically as $\Psi: \R^{N_\xi} \to \R$, where $N_\xi$ is the number of parameters required for parametrizing the positions and moments of the considered activation dipoles. The precise parametrization employed in our numerical tests is introduced in Section~\ref{sec:numerics}; at this stage, it is enough to note that the parametrization is smooth and such that no constraints are needed for the decision variables, meaning that one can resort to differentiation-based methods of global optimization. In particular, the needed first and second order derivatives can be calculated explicitly using the employed parametrization and applying basic matrix differentiation formulas to \eqref{eq:Aoptimal} and \eqref{eq:Doptimal}; see~\cite{Maaninen23} for more details.

\subsection{Minimization of the target functions}
Our standard algorithms for minimizing $\Psi$ are gradient descent and Newton's method accompanied by an inexact bisection line search that utilizes the Wolfe conditions. 
As the target function generally suffers from several local minima, the choice of the initial guess for the employed minimization algorithm plays a crucial role. The algorithm could,~e.g.,~be restarted from multiple initial guesses, and among the resulting designs, the one producing the smallest value for the optimization target could be chosen as the final minimizer. However, we simply resort to randomization or some heuristic for choosing a single initial guess in our numerical experiments.

The complete algorithm for optimizing an activation design is a combination of either Algorithm~\ref{alg:GD} (gradient descent) or Algorithm~\ref{alg:newton} (Newton's method) and Algorithm~\ref{alg:wolfe_alg}, which is an inexact bisection line search that is terminated when the Wolfe conditions are satisfied with predefined parameters. $H_\Psi$ denotes the Hessian of the optimization target and inequalities between (symmetric) matrices are interpreted in the sense of the partial ordering induced by positive-definiteness. In particular, the if-clause of Algorithm~\ref{alg:newton} guarantees that $H$ is strictly positive-definite, and $d$ is thus a descent direction for $\Psi$ at $\xi_i$.

\begin{algorithm}[t]
\caption{(Gradient descent)} 
\label{alg:GD}
\begin{algorithmic} 
\STATE Choose a tolerance $\epsilon > 0$, an initial guess $\xi_0$, the maximum number of iterations $N_{\rm GD}$, and a step size parameter $\lambda > 0$. Initialize $i = 0$. 
\WHILE {$|\nabla \Psi(\xi_i)| > \epsilon$ and $i < N_{\rm GD}$}
 \STATE \hspace{2mm} $\rhd$ Compute the search direction $d = \tilde{d}/|\tilde{d}|$, with $\tilde{d} = -\nabla \Psi(\xi_i)$. 
 \STATE \hspace{2mm} $\rhd$ Select the step size $\bar{\lambda}$ based on Algorithm~\ref{alg:wolfe_alg} with an initial guess $\lambda$.
 \STATE \hspace{2mm} $\rhd$ Set $\xi_{i+1} = \xi_i + \bar{\lambda}  d$.
   \STATE \hspace{2mm} $\rhd$ Update $i = i+1$.    
\ENDWHILE
\STATE {\bf return} $\xi_i$
\end{algorithmic}
\end{algorithm}

\begin{algorithm}[t]
\caption{(Newton's method)} 
\label{alg:newton}
\begin{algorithmic} 
\STATE Choose a tolerance $\epsilon > 0$, an initial guess $\xi_0$, the maximum number of iterations $N_{\rm Newton}$, a positiveness constraint $\delta > 0$, and a step size parameter $\lambda > 0$. Initialize $i = 0$. 
\WHILE {$|\nabla \Psi(\xi_i)| > \epsilon$ and $i < N_{\rm Newton}$}
\IF {$H_\Psi(\xi_i) \geq \delta I$ }
\STATE \hspace{2mm} $\rhd$ Set $H = H_\Psi(\xi_i)$
\ELSE
\STATE \hspace{2mm} $\rhd$ Set $H = H_\Psi(\xi_i) + (\delta - \mu)I$, where $\mu$ is the smallest eigenvalue of $H_\Psi(\xi_i)$.
 \ENDIF
 \STATE \hspace{2mm} $\rhd$ Compute the search direction $d = \tilde{d}/|\tilde{d}|$, with $\tilde{d} = -H^{-1}\nabla \Psi(\xi_i)$. 
 \STATE \hspace{2mm} $\rhd$ Select the step size $\bar{\lambda}$ based on Algorithm~\ref{alg:wolfe_alg} with an initial guess $\lambda$.
 \STATE \hspace{2mm} $\rhd$ Set $\xi_{i+1} = \xi_i + \bar{\lambda}  d$.
   \STATE \hspace{2mm} $\rhd$ Update $i = i+1$.    
\ENDWHILE
\STATE {\bf return} $\xi_i$
\end{algorithmic}
\end{algorithm}

\begin{algorithm}[t]
\caption{(Bisection method with the Wolfe conditions)}
\label{alg:wolfe_alg}
\begin{algorithmic}
\STATE Choose the maximum number of iterations $N_{\mathrm{Wolfe}}$ and a pair of scaling parameters $\beta_1, \beta_2 \in (0,1)$ satisfying $\beta_1 < \beta_2$. Initialize a counter as $l = 0$ and the lower and upper bounds for the bisection method with $\gamma_1 = 0$ and $\gamma_2 = \infty$. The initial step size 
$\lambda > 0$, the considered base point $\xi_i$, and the search direction $d$ are given as inputs.
\WHILE{$l < N_{\mathrm{Wolfe}}$}
    \IF{$\Psi(\xi_i + \lambda d) - \Psi(\xi_i) > \lambda \beta_1 \nabla \Psi(\xi_i)^\top d$}
        \STATE \hspace{2mm} $\rhd$ Set $\gamma_2 = \lambda$.
        \STATE \hspace{2mm} $\rhd$ Set $\lambda = \frac{1}{2}(\gamma_1 + \gamma_2)$.
    \ELSIF{$\nabla \Psi(\xi_i + \lambda d)^\top d < \beta_2 \nabla \Psi(\xi_i)^\top d$}
        \STATE \hspace{2mm} $\rhd$ Set $\gamma_1 = \lambda$.
        \IF{$\gamma_2 = \infty$}
            \STATE \hspace{2mm} $\rhd$ Set $\lambda = 2\gamma_1$.
        \ELSE
            \STATE \hspace{2mm} $\rhd$ Set $\lambda = \frac{1}{2}(\gamma_1 + \gamma_2)$.
        \ENDIF
    \ELSE
        \RETURN $\bar{\lambda} = \lambda$
    \ENDIF
    \STATE \hspace{2mm} $\rhd$ Set $l = l + 1$.
\ENDWHILE
\RETURN $\bar{\lambda}  = \lambda$
\end{algorithmic}
\end{algorithm}

The parameter $\beta_1$, which is used in the first Wolfe condition, describes a sufficient decrease in the value of the target function: the step size is accepted if it results in a decrease in the target function that is at least $\beta_1$  times the decrease predicted by the directional derivative at the base point of the line search. The parameter $\beta_1$ is often chosen to be quite small, which is already enough for guaranteeing a decrease in the target function. The parameter $\beta_2$ controls the so-called curvature condition. If the slope at the point proposed by the step size is larger than $\beta_2$ times the initial slope, the curvature condition is satisfied. As the initial slope is guaranteed to be negative by Algorithms~\ref{alg:GD} and \ref{alg:newton}, this means that the rate of decrease at the proposed point can be at most $\beta_2$ times the rate of decrease at the base point of the line search. The parameter $\beta_2$ is often chosen to be quite close to one.

\section{Sequential edge-promoting optimization of activations}
\label{sec:opti}

If a Gaussian prior,~i.e.,~case 1 in Section~\ref{ssec:A_optimality} is considered, all essential material for implementing our algorithm for optimizing the activation designs has already been included in Section~\ref{sec:oed}. However, if case 2 of Section~\ref{ssec:A_optimality} is tackled, the optimization routine must at each iteration be combined with a lagged diffusivity iteration for computing an approximation of the MAP estimate for the posterior \eqref{eq:post_u} with $\Phi$ defined by the smoothened TV regularizer \eqref{eq:aRRa}. The purpose of this section is to summarize this combined procedure as a concise algorithm,~i.e.~Algorithm~\ref{alg:basic_optimization}, which is essentially the same as the one presented in \cite[Section~5]{Helin22}.

\begin{algorithm}[t]
\caption{(Edge-promoting sequential design)}
\label{alg:basic_optimization}
\begin{algorithmic}
  \STATE{Select the prior parameters $T>0$ and $\gamma >0$ and a tolerance for the stopping criterion $\tau > 0$.
    }

\vspace{1.5mm}

  \STATE{Initialization:}
  \STATE{\hspace{2mm} $\rhd$ Set $\widehat{c}_0 = \mathbf{1} \in \R^{N_c}$.}
  \STATE{\hspace{2mm} $\rhd$ Define $\Gamma_{0} := \Theta(\widehat{c}_0)^{-1}$ according to \eqref{eq:H}.}

\vspace{1.5mm}

  \STATE{Iteration:}

  \FOR{$k=1, \dots, N_{a}$}
  \STATE{\hspace{2mm} $\rhd$ Solve for the next activation $p_{k} := p(\xi_k)$, where $\xi_k$ is a solution to \eqref{eq:Aoptimal} or \eqref{eq:Doptimal}, with $\Sigma(\xi)$ defined via \eqref{eq:post_seq}, deduced by Algorithm~\ref{alg:GD}, Algorithm~\ref{alg:newton} or an exhaustive search (cf.~\cite{Burger21}).}
  \STATE{\hspace{2mm} $\rhd$ Form the system matrix $K(p_k)$ and `measure' the data $y_k$.}
  \STATE{\hspace{2mm} $\rhd$ Set $j=0$, $\widehat{c}_{k}^{(0)}=\widehat{c}_{k-1}$, and $\Delta \Phi = \tau + 1$.}
  \WHILE{$\Delta \Phi > \tau$}
  \STATE{\hspace{2mm} $\rhd$ Set $j = j+1$.}
     \STATE{\hspace{2mm} $\rhd$ Compute $\widehat{c}_{k}^{(j)}$ according to \eqref{eq:step2}.}
     \STATE{\hspace{2mm} $\rhd$ Compute $\Delta \Phi = |\Phi(\widehat{c}_{k}^{(j-1)}) - \Phi(\widehat{c}_{k}^{(j)})|/\Phi(\widehat{c}_{k}^{(j)})$.}
  \ENDWHILE
  \STATE{\hspace{2mm} $\rhd$ Define $\widehat{c}_k = \widehat{c}_{k}^{(j)}$ as well as $\Gamma_{k} = \Gamma_{k}^{(j)}$ formed according to \eqref{eq:posterior_cov}}.
  \ENDFOR

  \vspace{2mm}

  \RETURN $\widehat{c}_{N_a}$, $\Gamma_{N_a}$, and $(p_1, \dots, p_{N_a})$.
\end{algorithmic}
\end{algorithm}

In addition to the sequentially optimized activation $p_1, \dots, p_{N_a}$, Algorithm~\ref{alg:basic_optimization} also returns the final reconstruction $\widehat{c}_{N_a}$ and the associated spread estimator $\Gamma_{N_a}$,~i.e.,~the mean and the covariance of the final approximative Gaussian posterior for the MNP concentration after $N_a$ activations. A motivation for the stopping criterion of the lagged diffusivity iteration in the inner loop can be found in~\cite{Arridge13,Helin22}. Note that we also employ an exhaustive algorithm for computing the optimal sequential designs in Algorithm~\ref{alg:basic_optimization} to be able to deduce upper limits for the performance of this adaptive edge-promoting approach to the activation design in MRXI. We refer to \cite{Burger21,Helin22} for considerations on efficient implementation of such an exhaustive optimization routine in the framework X-ray imaging.

\section{Numerical experiments}
\label{sec:numerics}
All numerical experiments are performed in a two-dimensional setting, with $\Omega \subset \R^2$ being a disk of radius $\rho>0$ centered at the origin. In particular, all spatial variables, dipole moments and sensor orientations introduced in Section~\ref{sec:mm} are modeled as two-dimensional vectors. This can be considered an approximation for a (somewhat unrealistic) setting, where the MNPs are concentrated around a single cross-section of a three-dimensional body, modeled by $\Omega$, and the sensors and activations lie in the two-dimensional plane defined by $\Omega$, with their orientations and dipole moments being parallel to that same plane. 

To be more precise, the sensors and activations are placed on a circle with radius $1.1 \rho$, and their orientations and activations are modeled by unit vectors. In all tests, there are $N_s = 36$ sensors that are equiangularly spaced and oriented toward the origin. The position and moment of an activation are parametrized by two angular variables, with one of them defining the position on the measurement circle and the other giving the orientation of the activation. In particular, the decision variable $\xi$ considered in Section~\ref{sec:oed} is a concatenation of two or several such angular variables, whose periodicity enables using gradient descent and Newton's method for estimating optimal designs without any additional constraints. However, observe that reversing the orientation of an activation only changes the signs of all associated (noiseless) measurements, which does not affect the A- or D-optimality of the design in question. Hence, the deduced optimal designs are post-processed for visualization purposes by reversing the orientations of those activations that originally point toward the interior of the measurement circle.

The integral kernel in \eqref{eq:complete_nondisc_measurement} depends explicitly and smoothly on the parameter pair $(a, \alpha)$ defining the position and moment of an activation dipole, and thus the kernel also depends explicitly and smoothly on the two angular decision variables defining an activation in our parametrization. In particular, one can straightforwardly, yet tediously deduce the first and second derivatives of the elements of the measurement matrix $\mathbf{K}(\mathbf{p}_k)$ in \eqref{eq:ax_final} with respect to the angular decision variables defining the components of $\mathbf{p}_k$. The gradient and Hessian of the optimization target, defined by \eqref{eq:Aoptimal} or \eqref{eq:Doptimal}, can thus also be given explicitly by resorting to basic rules of matrix differentiation, independently of whether one or several activations are optimized simultaneously. We do not write down the formulas for these derivatives here but refer to \cite{Maaninen23} for more details.

The measurement noise is assumed to have independent components with a common standard deviation $\eta > 0$. The disk $\Omega$ is discretized by initially constructing a uniform pixel grid for an origin-centered square with side length $2\rho$ and then removing the pixels that lie outside $\Omega$. The MNP concentration is modeled as a piecewise constant function with respect to such a pixelification of $\Omega$, with $c$ carrying the corresponding constant concentration values. Three levels of discretization are used: the number of components in $c$ is $N_c = 4293$ for computing reconstructions, a coarser grid with $N_{\rm opt} = 901$ pixels is used for optimizing the designs, while measurement data is simulated with $N_{\rm data} = 5140$ degrees of freedom to avoid an inverse crime. We refer to \cite{Burger21,Helin22} for technical details on performing the optimization of the measurement design on a sparser grid by resorting to interpolation between pixelifications.

The free parameters in Algorithms~\ref{alg:GD}, \ref{alg:newton}, \ref{alg:wolfe_alg} and \ref{alg:basic_optimization} are set to $\epsilon = 10^{-5}$, $\lambda = 1$, $N_{\rm GD} = N_{\rm Newton} = 50$, $\delta = 10^{-5}$, $N_{\rm Wolfe} = 15$, $\beta_1 = 10^{-10}$, $\beta_2 = 0.9$, $T=10^{-6}$, $\gamma = 10$ and $\tau = 10^{-3}$. The effect of the last three parameters on the performance of Algorithm~\ref{alg:basic_optimization} is discussed in \cite{Helin22}, where the same values are used in X-ray imaging. The low value for $\beta_1$ means that even a small decrease in the objective function is considered sufficient in the bisection line search. Due to the existence of several local minima for the A- and D-optimality objective functions (cf.~\cite{Burger21,Maaninen23}), the parameters controlling Algorithms~\ref{alg:GD}, \ref{alg:newton} and \ref{alg:wolfe_alg} may have a considerable effect on the proposed experimental designs since gradient descent and Newton's method can,~e.g.,~converge to different local minima even when starting from the same initial guess. We do not claim that the aforelisted parameter choices are optimal. If one resorts to exhaustive search in Algorithm~\ref{alg:basic_optimization}, the value of the optimization target is evaluated on a equidistant grid of $100 \times 100$ points over $[0, 2 \pi)\times [0, 2 \pi)$ to deduce the decision variable pair that minimizes the considered optimization target.

\begin{remark}
Although we have ignored many physical parameters in Section~\ref{sec:mm} and the orientations/moments of the sensors/activations have been scaled to unit length, there still remain two parameters that can be tuned in the experiments: the radius $\rho$ of $\Omega$ and the noise level $\eta$. As the ignored physical parameters and the magnitudes of the activations only scale the measurements, controlling the domain size and the noise level still provide sufficient flexibility for modeling the physical properties of the measurement setup.
\end{remark}

\subsection{Gaussian prior}
Our first two experiments consider a Gaussian prior and simultaneous optimization of $N_a = 10$ activations. The prior covariance is given elementwise as 
\begin{equation}
 \label{eq:prior_cov}
 (\Gamma_{0})_{i,j} = \gamma^2 \exp \left(-\frac{| x_i - x_j |^2}{2\ell^2} \right),
 \end{equation}
where $\ell>0$ is the so-called correlation length, $\gamma>0$ is the pixelwise standard deviation, and $x_i$ and $x_j$ are the coordinates of the pixels with indices $i$ and $j$. We employ the same standard deviation $\gamma = 1$ in all experiments with a Gaussian prior, but the correlation length varies between individual tests.

\subsubsection{Gaussian Test~1: a homogeneous disk}
In the first experiment, we select the standard deviation of the noise, the radius of the imaged domain and the correlation length for the prior as $\eta = 1$, $\rho = 0.5$ and $\ell = 0.15$, respectively. The initial positions and orientations of the activations are independently drawn from the uniform density over $[0, 2\pi]$, and Algorithm~\ref{alg:newton},~i.e.,~Newton's method, is used for finding both A- and D-optimal designs. For the A-optimality target function defined in \eqref{eq:Aoptimal}, the weight $A$ is chosen to be the identity matrix, meaning that the reconstruction quality is considered equally important everywhere in $\Omega$. 

\begin{figure}
    \centering
    \includegraphics[width = 0.3\columnwidth]{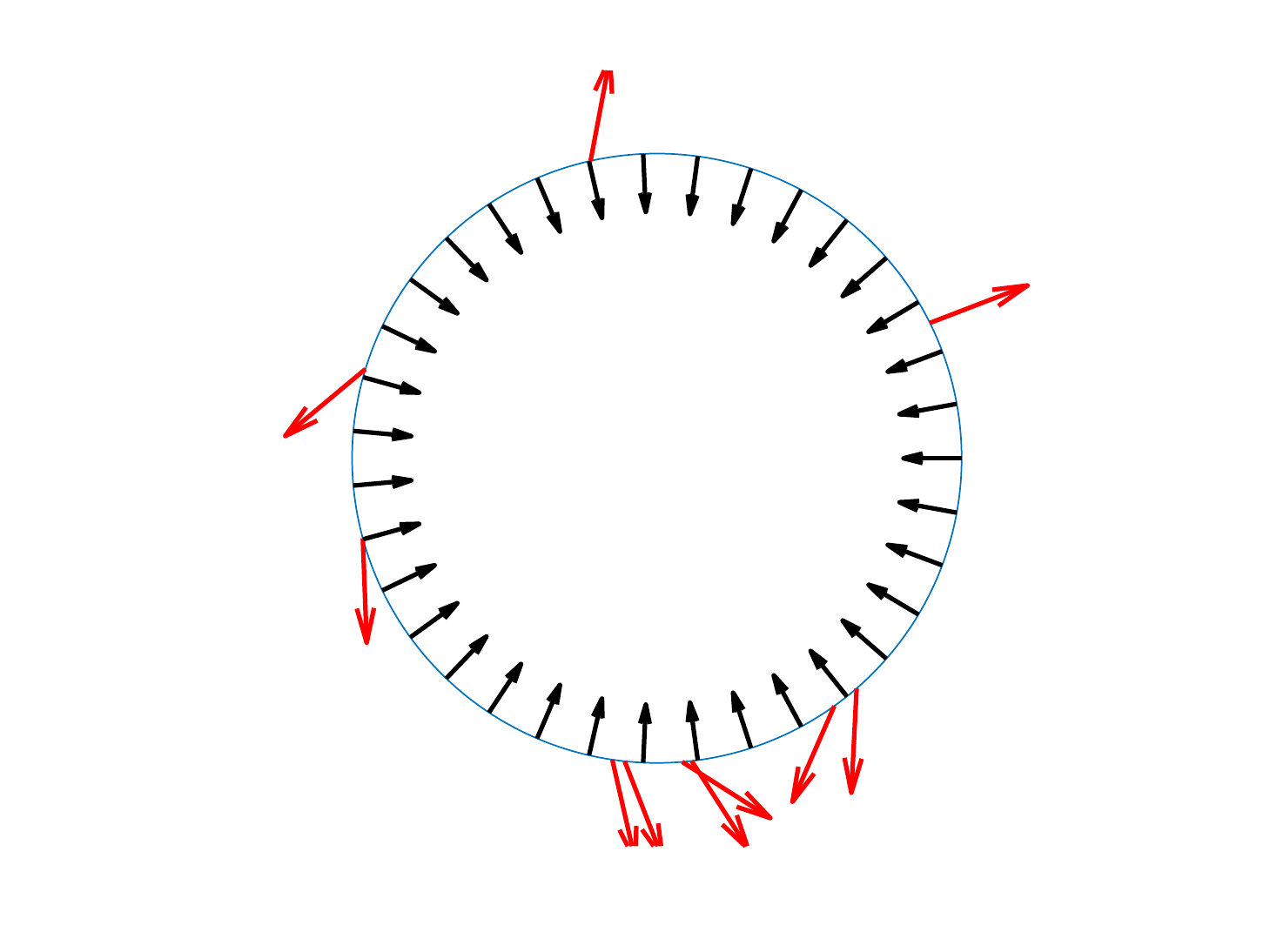}
    \includegraphics[width = 0.3\columnwidth]{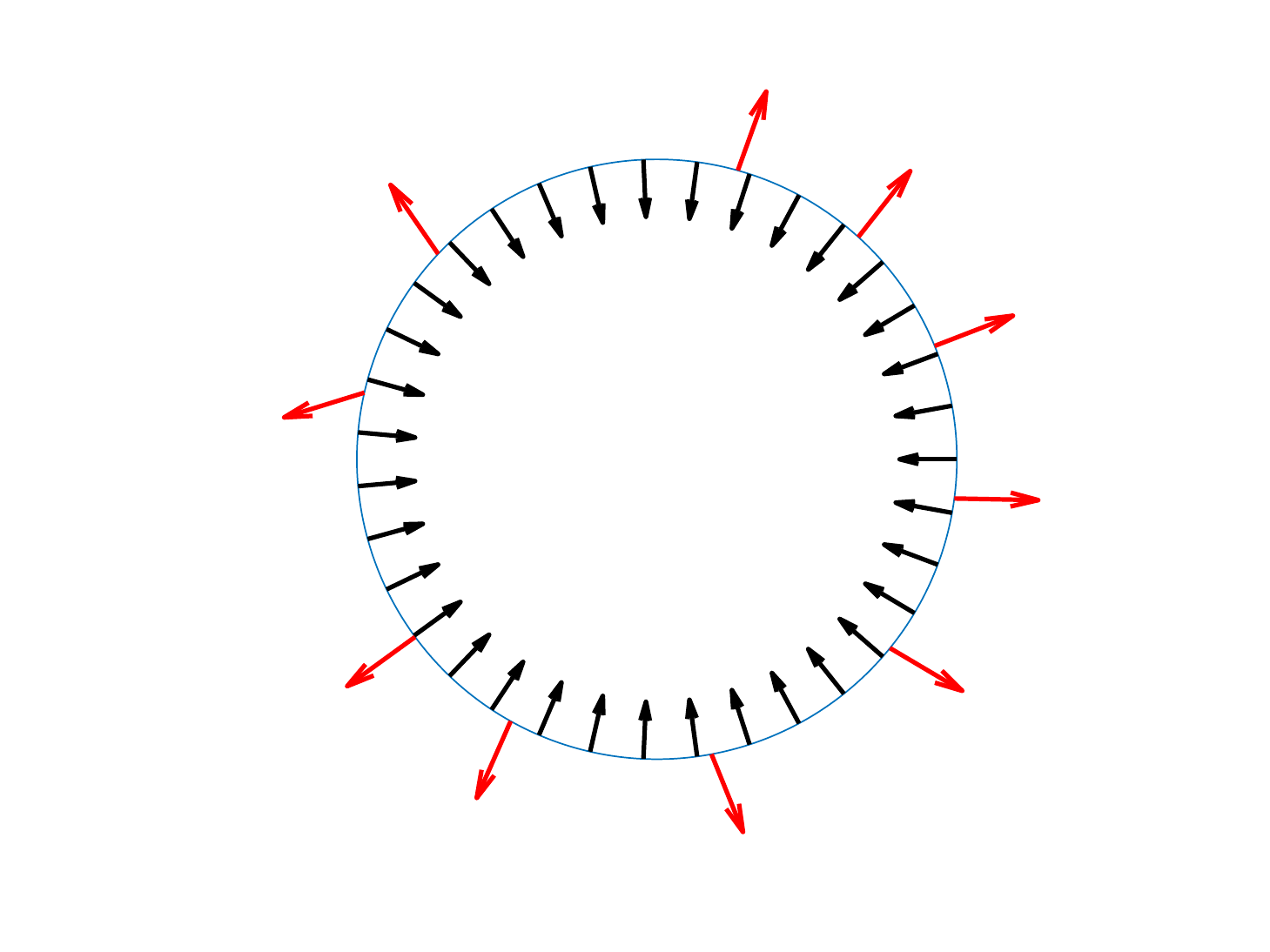}
    \includegraphics[width = 0.3\columnwidth]{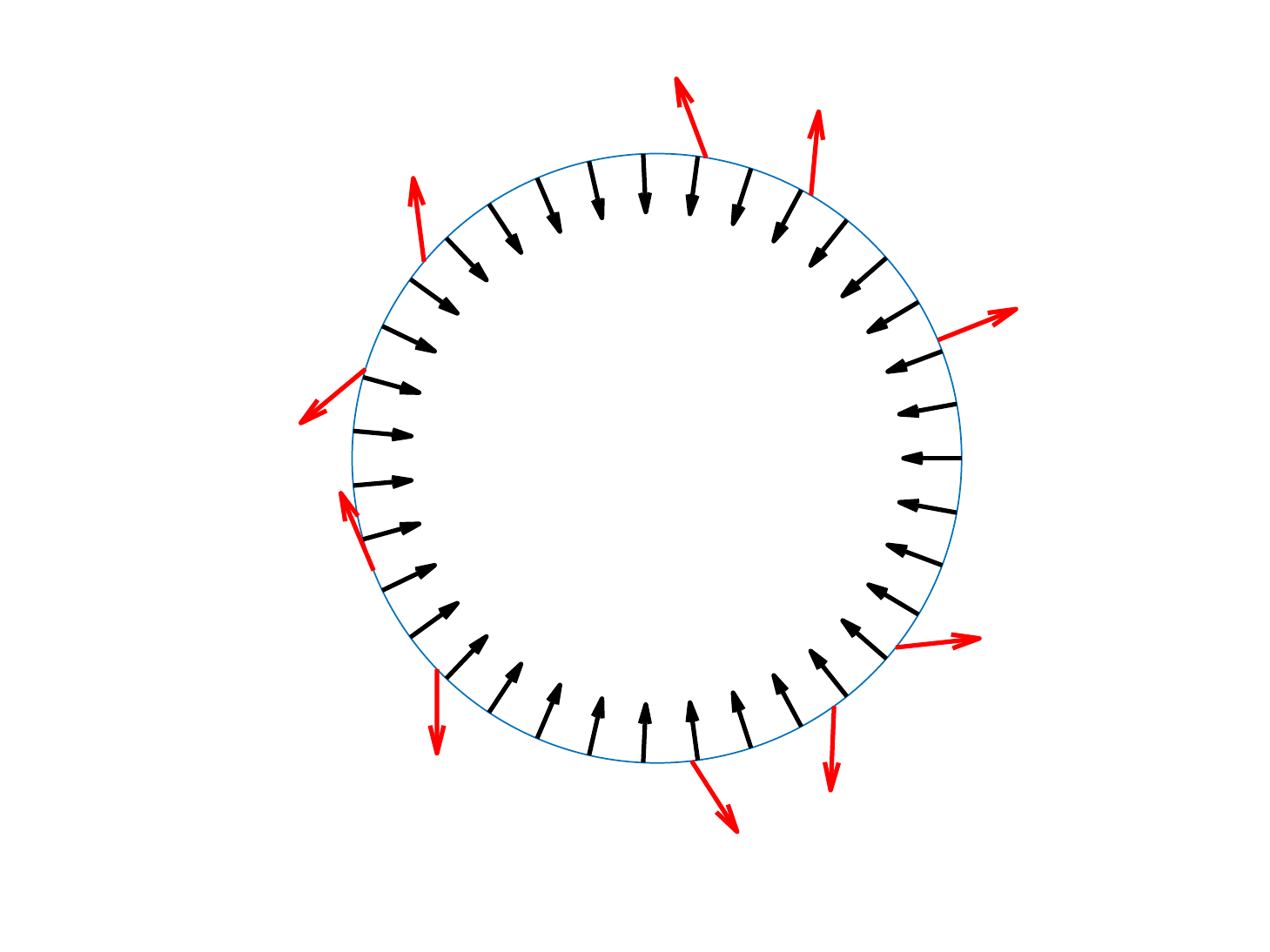}
    \caption{{\sc Gaussian~Test~1.} Left: initial activation configuration. Center: approximate A-optimal design. Right: approximate D-optimal design. The black arrows depict the sensors and the red arrows the optimized activations.} 
    \label{fig:gauss_test1_1}
\end{figure}

The initial activation configuration is visualized in the left-hand image of Figure~\ref{fig:gauss_test1_1}, and the middle and right images show the estimated A- and D-optimal activation designs, respectively. For both A- and D-optimality, the optimized activations lie approximately equidistantly around the object and point roughly away from its center (after a possible reorientation by 180 degrees), though in the D-optimal design the directions vary considerably more. This suggests that in the case of nonadaptivity and a radially symmetric target, placing the activations symmetrically is a reasonable approach, which cannot be considered very surprising.

\begin{figure}
    \centering
    \includegraphics[width = 0.49\columnwidth]{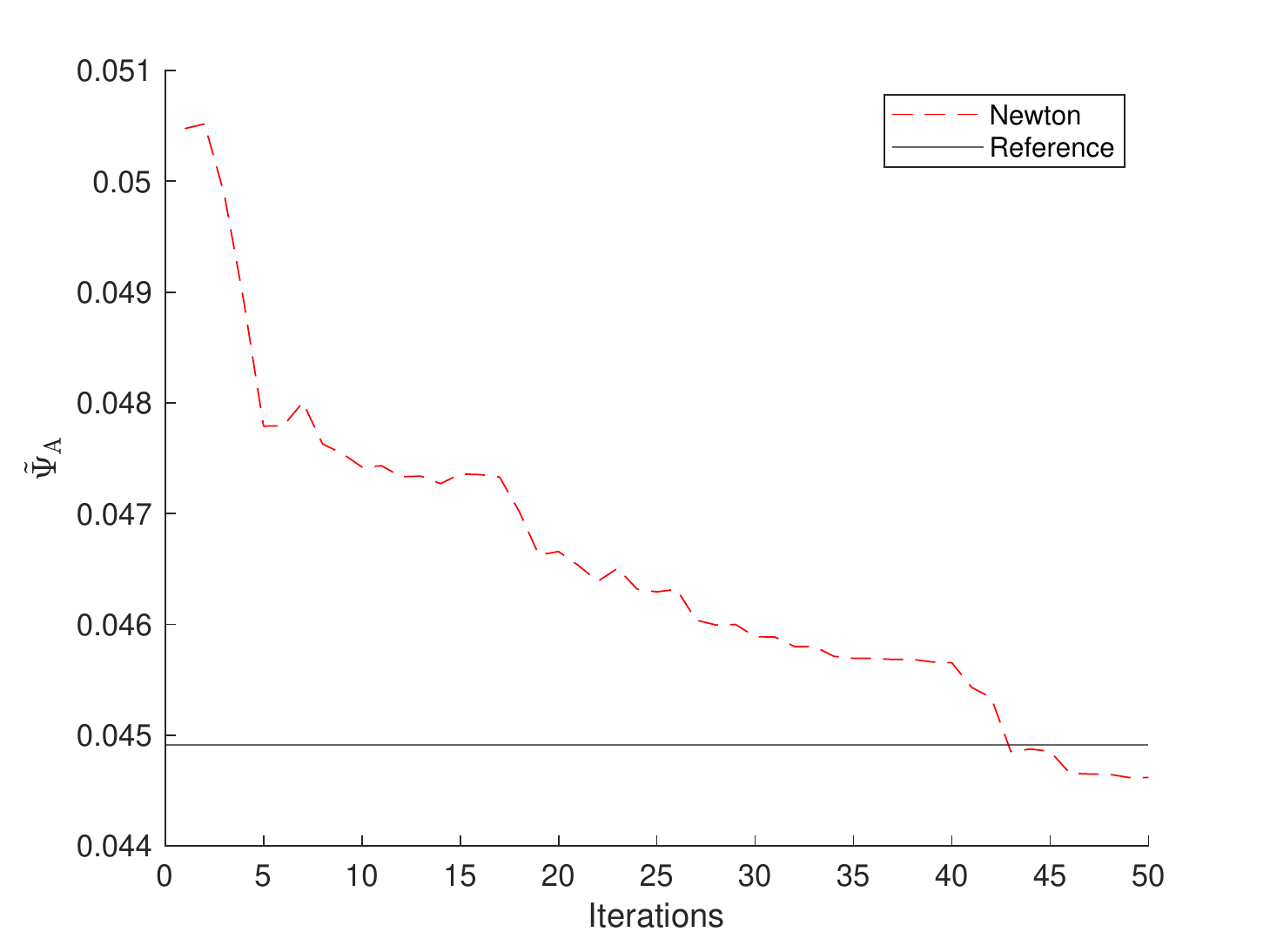}
    \includegraphics[width = 0.49\columnwidth]{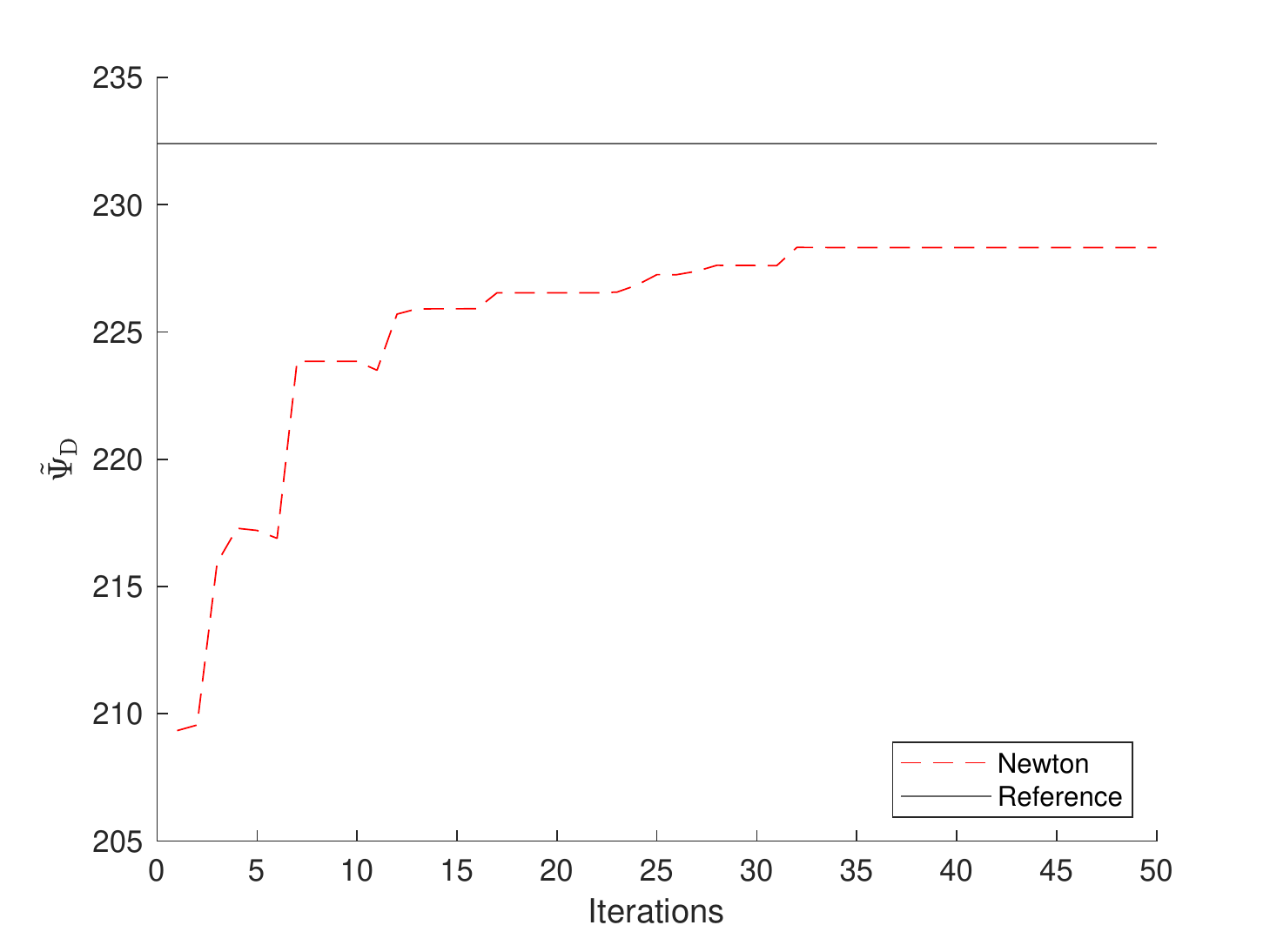}
    \caption{{\sc Gaussian~Test~1.} Left: the expected $L^2(\Omega)$ reconstruction error $\tilde{\Psi}_{\rm A}$ as a function of the iteration number when deducing the A-optimal design. Right: the information gain $\tilde{\Psi}_{\rm D}$ as a function of the iteration number when deducing the D-optimal design. The horizontal lines depict the values of $\tilde{\Psi}_{\rm A}$ (left) and $\tilde{\Psi}_{\rm D}$ (right) for equidistant activations pointing away from the center of $\Omega$.} 
    \label{fig:gauss_test1_2}
\end{figure}

The left-hand image in Figure~\ref{fig:gauss_test1_2} shows the expected $L^2(\Omega)$ reconstruction error $\tilde{\Psi}_{\rm A}$, defined in \eqref{eq:L2_error}, as a function of the iteration number of Newton's method. Similarly, the right-hand image illustrates the evolution of the properly scaled information gain (cf.~\cite{Burger21})
$$
\tilde{\Psi}_{\rm D}(\xi) = \frac{1}{2}\big( \log ( \det \Gamma_0) - \tilde{\Psi}_{\rm D}(\xi) \big)
$$
during the optimization procedure for a D-optimal design. The horizontal lines in Figure~\ref{fig:gauss_test1_1} depict the values of $\tilde{\Psi}_{\rm A}$ and $\tilde{\Psi}_{\rm D}$ for a precisely equidistant configuration of activations pointing exactly away from the center of $\Omega$. For both A- and D-optimality, the optimized target function values are close to those for the reference design, with the deduced A-optimal design even slightly surpassing the reference value, possibly thanks to the slight variations in the activation directions. The final D-optimal design corresponds to a lower information gain than the reference design, which suggests the optimization process got stuck in a local optimum.

\subsubsection{Gaussian Test~2: a semidiscoidal ROI}
\begin{figure}
    \centering
    \begin{tikzpicture}[scale=1.2]
      \useasboundingbox (-1.5,-1.6) rectangle (1.5,1.5);
      \draw[fill=green!20!white] (0, 0) circle (1) node[left=3mm] {ROI};
      \draw[fill=red!20!white] (0,-1) arc(-90:90:1) -- cycle;
      \draw[thick] (0.0, -1.0) -- (0.0, 1.0);
    \end{tikzpicture}
    \includegraphics[width = 0.35\columnwidth]{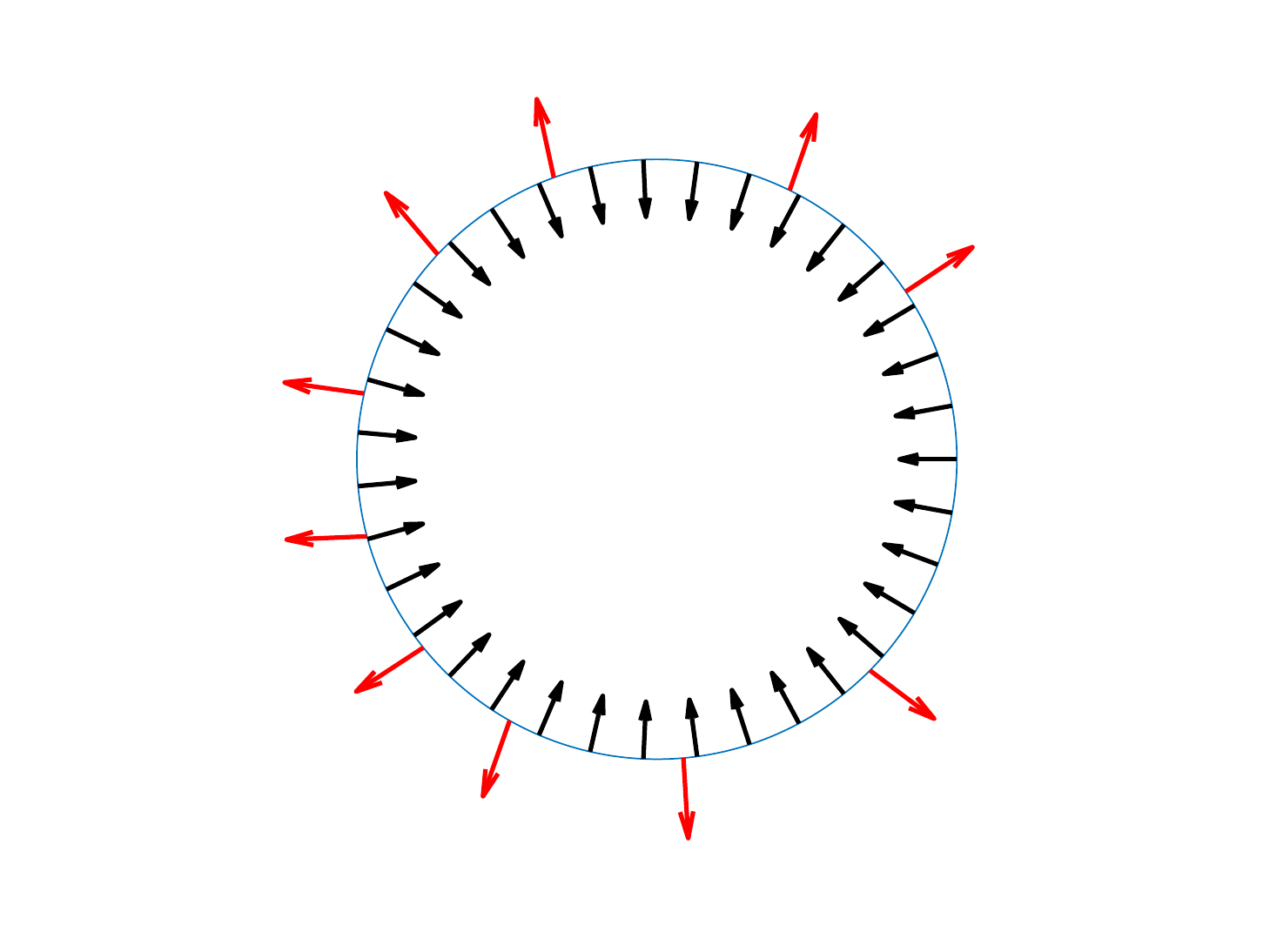}
    \includegraphics[width = 0.35\columnwidth]{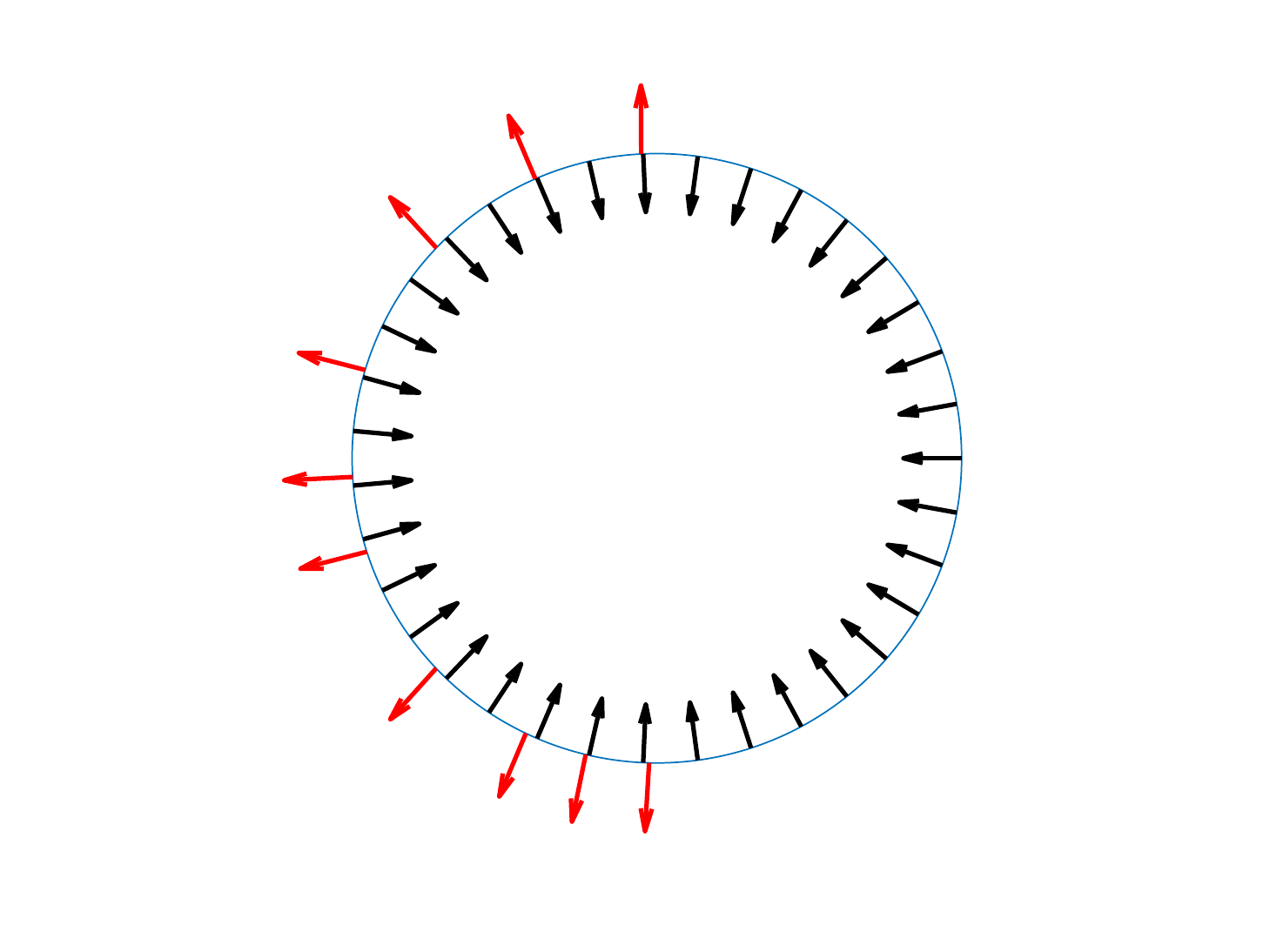}
    \caption{{\sc Gaussian~Test~2.} Left: ROI is the left half of $\Omega$. Center: approximate A-optimal design for $\eta = 1$, $\rho = 0.5$ and $\ell = 0.15$. Right: approximate A-optimal design for $\eta = 0.1$, $\rho = 5$ and $\ell = 1.5$. The black arrows depict the sensors and the red arrows the optimized activations. The images are not in scale, but the diameter of the right-hand version of $\Omega$ is ten times larger than of that at the middle.}
    \label{fig:gauss_test2_1}
\end{figure}
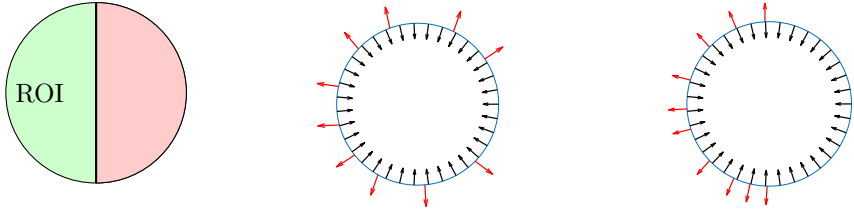

Our second test with a Gaussian prior investigates a setting where the ROI is the left half of $\Omega$, as illustrated in Figure~\ref{fig:gauss_test2_1}. Only A-optimality is considered. The initial guess for the set of activations is the same random configuration as used in the first test, shown on left in Figure~\ref{fig:gauss_test1_1}. Newton's method is employed for finding approximate A-optimal configurations for two sets of parameter values: $\eta = 1$, $\rho = 0.5$ and $\ell = 0.15$ as in the first test, and $\eta = 0.1$, $\rho = 5$ and $\ell = 1.5$. The resulting A-optimal activation configurations are visualized in the middle and right-hand images of Figure~\ref{fig:gauss_test2_1}, respectively. In both cases, the activations are again oriented approximately away from the center of the domain, but their locations have clearly shifted toward the ROI on the left-hand side of $\Omega$. This phenomenon is more pronounced for the latter parameter combination, presumably due to the larger size of the imaged object that makes activations on the side opposite to the ROI less useful.

\begin{figure}
    \centering
   \includegraphics[width = 0.49\columnwidth]{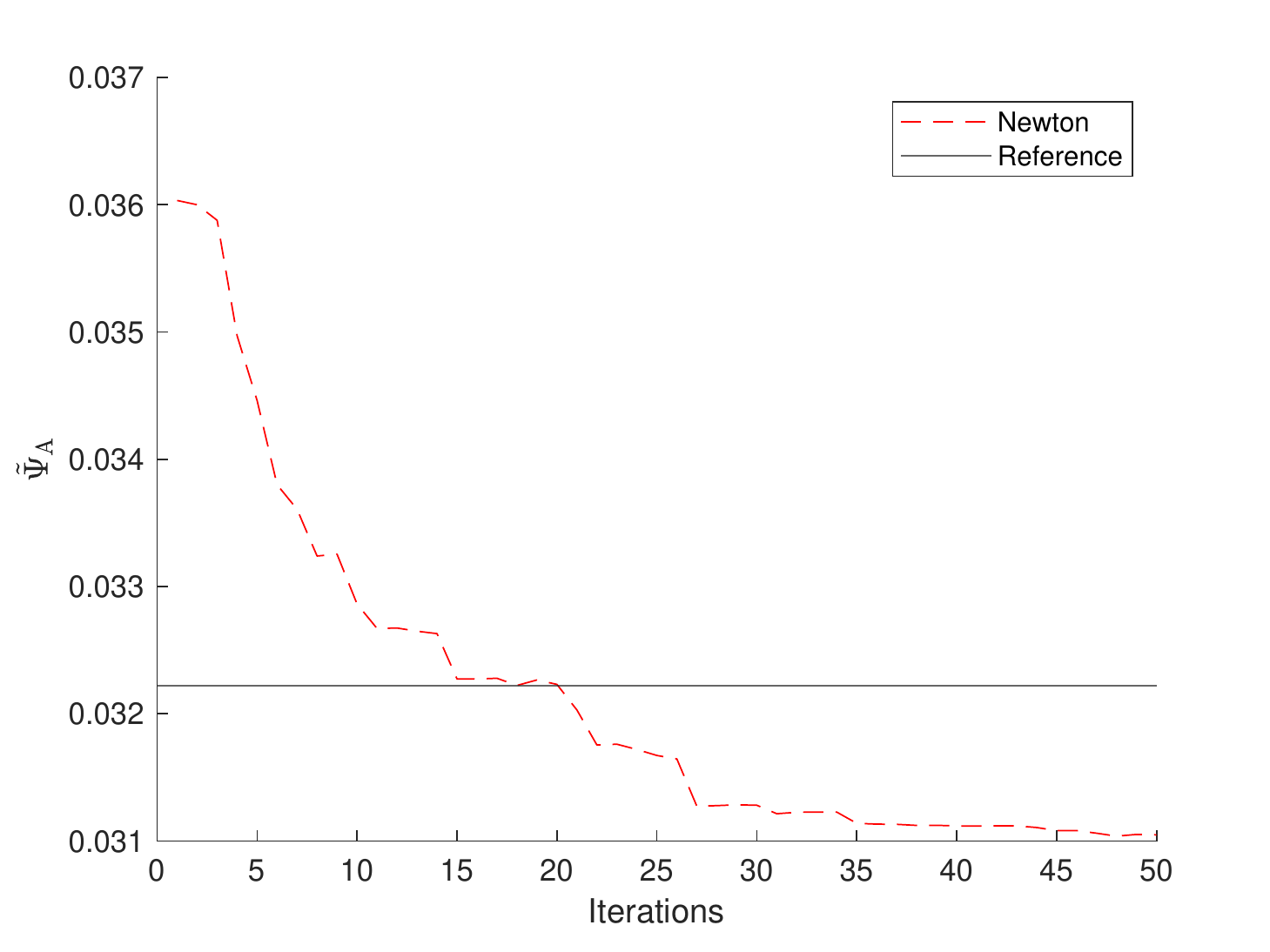}
    \includegraphics[width = 0.49\columnwidth]{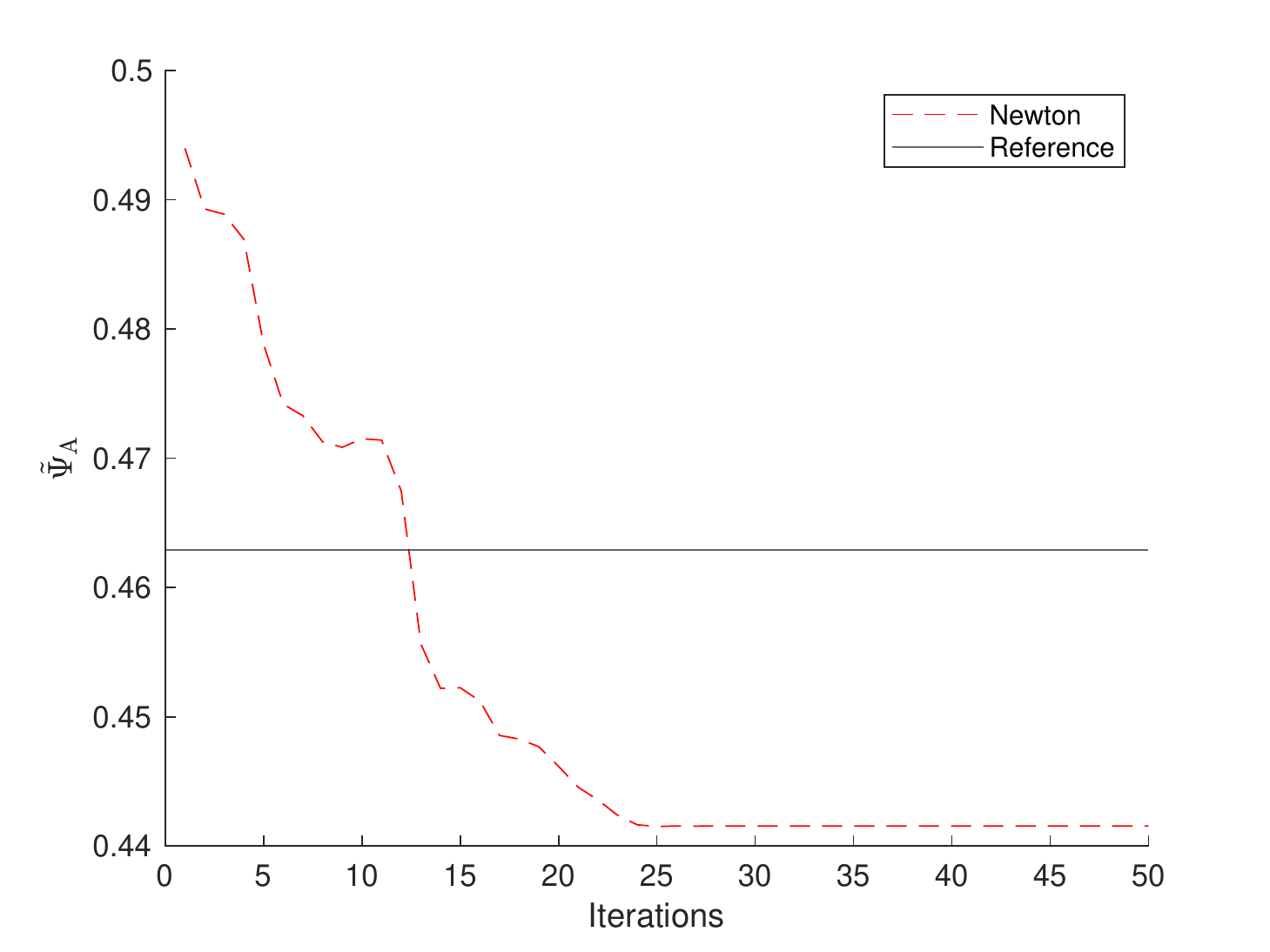}
    \caption{{\sc Gaussian~Test~2.} The expected $L^2(\Omega)$ reconstruction error $\tilde{\Psi}_{\rm A}$ as a function of the iteration number. Left: $\eta = 1$, $\rho = 0.5$ and $\ell = 0.15$. Right: $\eta = 0.1$, $\rho = 5$ and $\ell = 1.5$. The horizontal lines depict the values of $\tilde{\Psi}_{\rm A}$ for the considered set of parameters and equidistant activations pointing away from the center of $\Omega$.} 
    \label{fig:gauss_test2_2}
\end{figure}

Figure~\ref{fig:gauss_test2_2} illustrates the evolution of the expected $L^2(\Omega)$ reconstruction error $\tilde{\Psi}_{\rm A}$ during Newton's iterations for the two sets of parameter values. Again, the horizontal lines denote the expected errors for a reference design consisting of equally spaced activations oriented away from the center of $\Omega$. The optimized configurations show a greater advantage over the reference design when compared to the experiment with the ROI covering all of $\Omega$. This is to be expected since the reference design has not been adapted to the current setup where only the left side of the object is of interest.

\subsection{TV prior and sequential designs}
\begin{figure}
	\centering
  \includegraphics[width = 0.39\columnwidth]{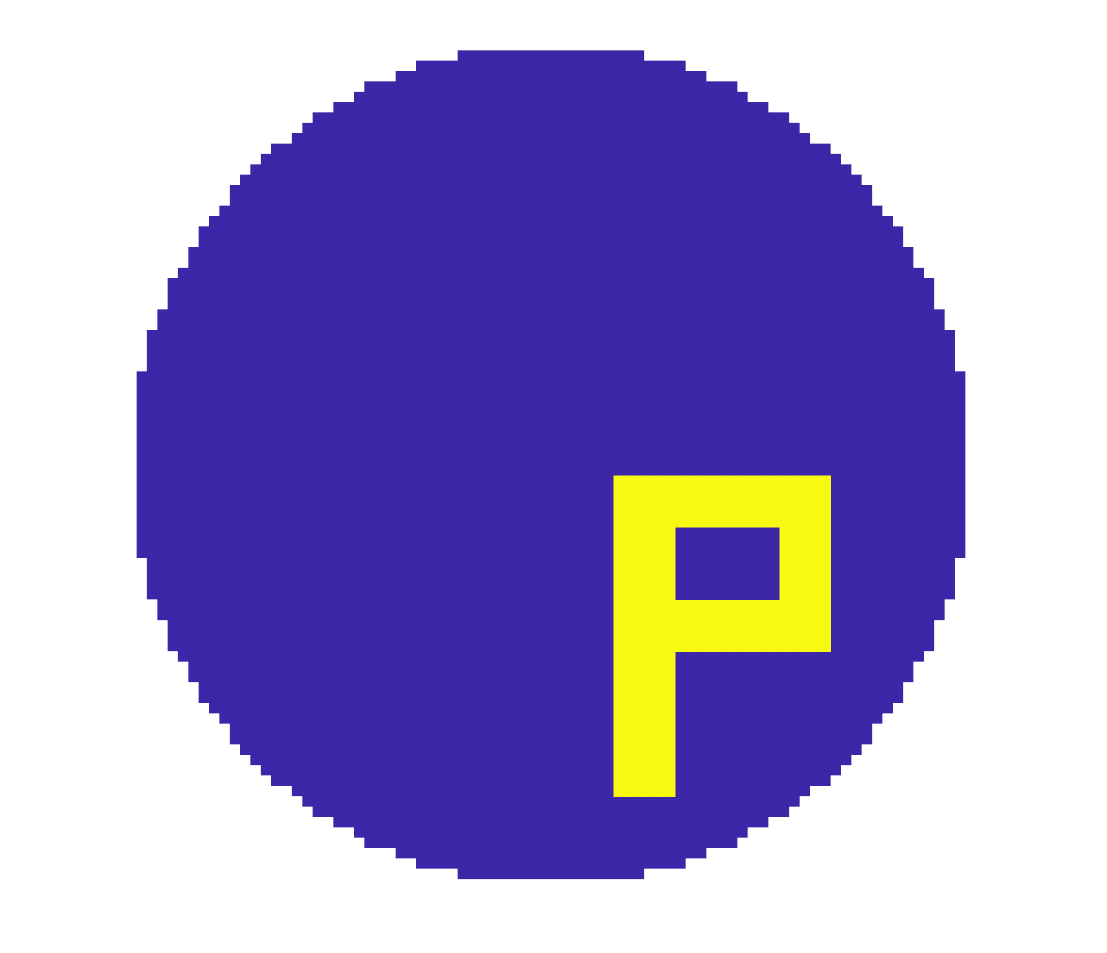}
  \includegraphics[width = 0.59\columnwidth]{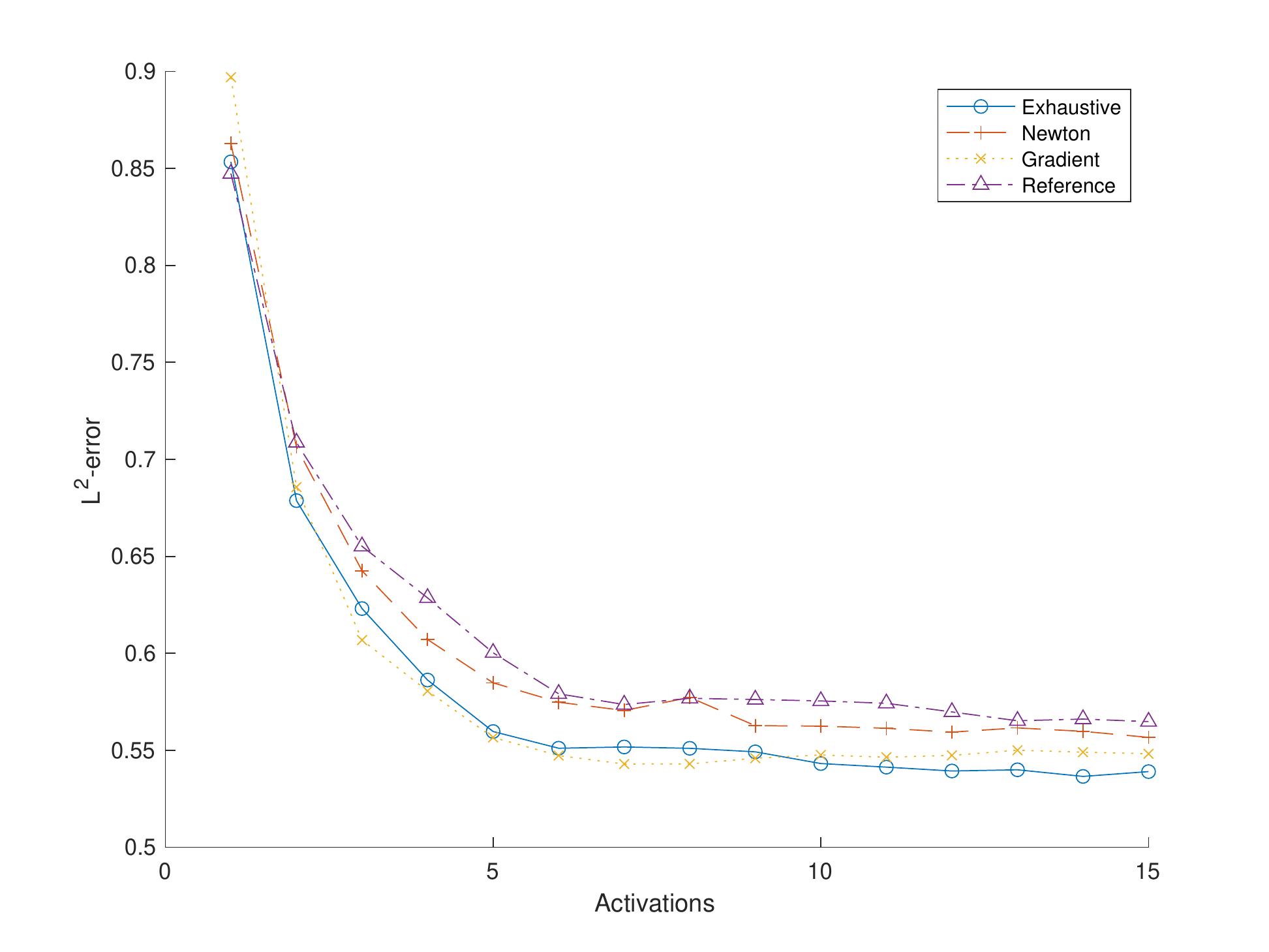}
	\caption{{\sc TV~Test~1.} Left: target MNP concentration. Right: relative $L^2(\Omega)$ reconstruction errors as functions of the number of activations in Algorithm~\ref{alg:basic_optimization}  when the specifications of the activations are deduced by Newton's method, Gradient descent, exhaustive search and the deterministic reference procedure.}
	\label{fig:TV_test1_1}
\end{figure}

Next, we apply the adaptive edge-promoting Algorithm~\ref{alg:basic_optimization} to the piecewise constant test phantom illustrated in the left-hand image of Figure~\ref{fig:TV_test1_1}, with two different choices for the noise-size parameter pair $(\eta, \rho)$. Only A-optimality is considered in this adaptive sequential context.

\subsubsection{TV Test~1: P-shaped inclusion}
The target is characterized by a homogeneous background contaminated by a P-shaped inhomogeneity located in the bottom-right quadrant of $\Omega$ and carrying unit MNP density. The size of the domain and the measurement noise level are defined by the parameter values and $\rho = 0.5$ and $\eta = 0.1$, respectively. The target is deliberately chosen to demonstrate advantages of Bayesian OED: if the P-shaped inclusion were positioned close to the center of $\Omega$, the optimal design would most likely approximately correspond to activations placed equiangularly around the object, and the optimization would not provide a clear benefit. However, considering a target that exhibits interesting features only in a small subregion of $\Omega$ close to its boundary leads to optimal designs where the activations are concentrated close to the particular area of interest.

\begin{figure}
	\centering
  \includegraphics[width = 0.3\columnwidth]{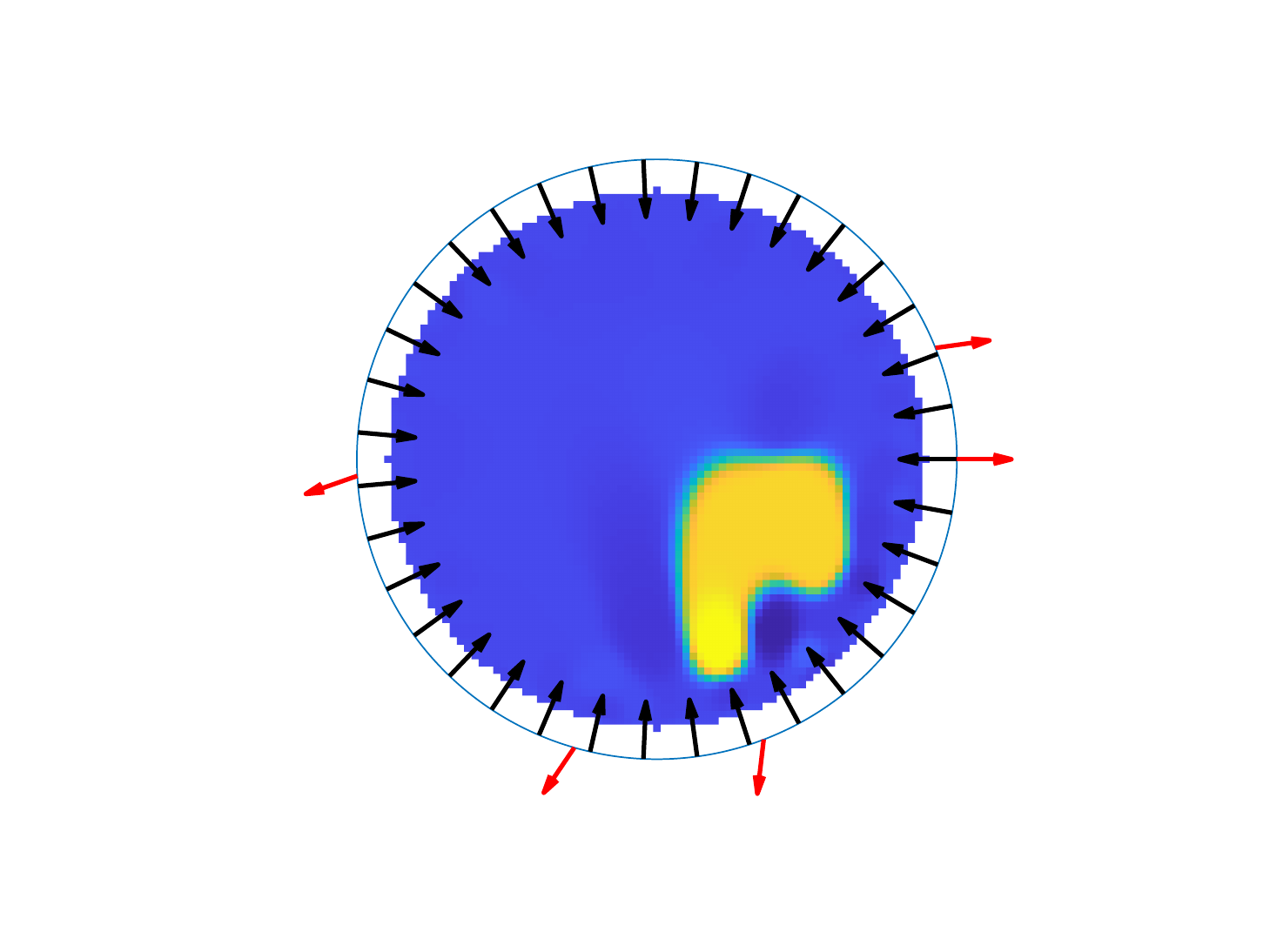}
  \includegraphics[width = 0.3\columnwidth]{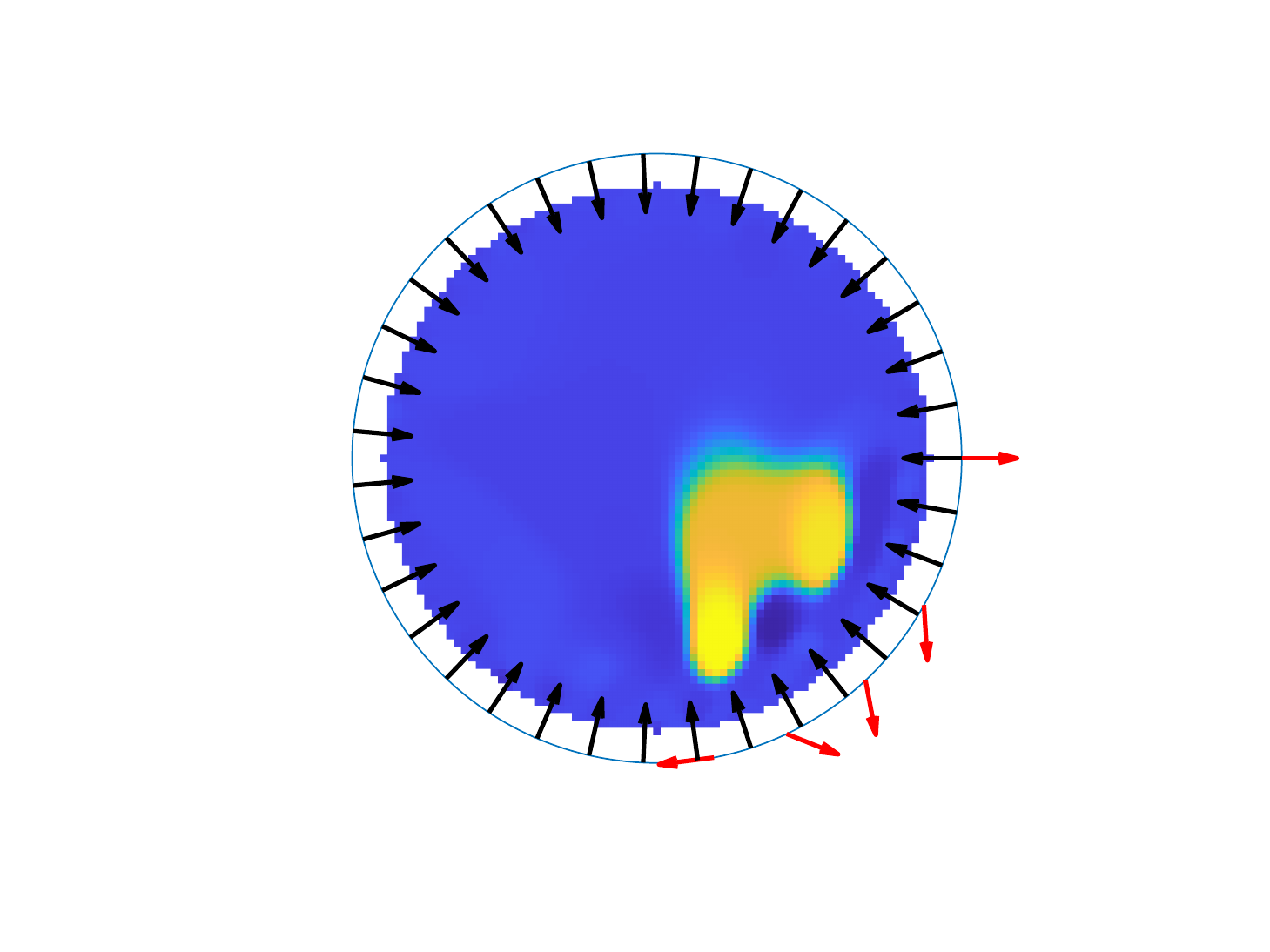}
  \includegraphics[width = 0.3\columnwidth]{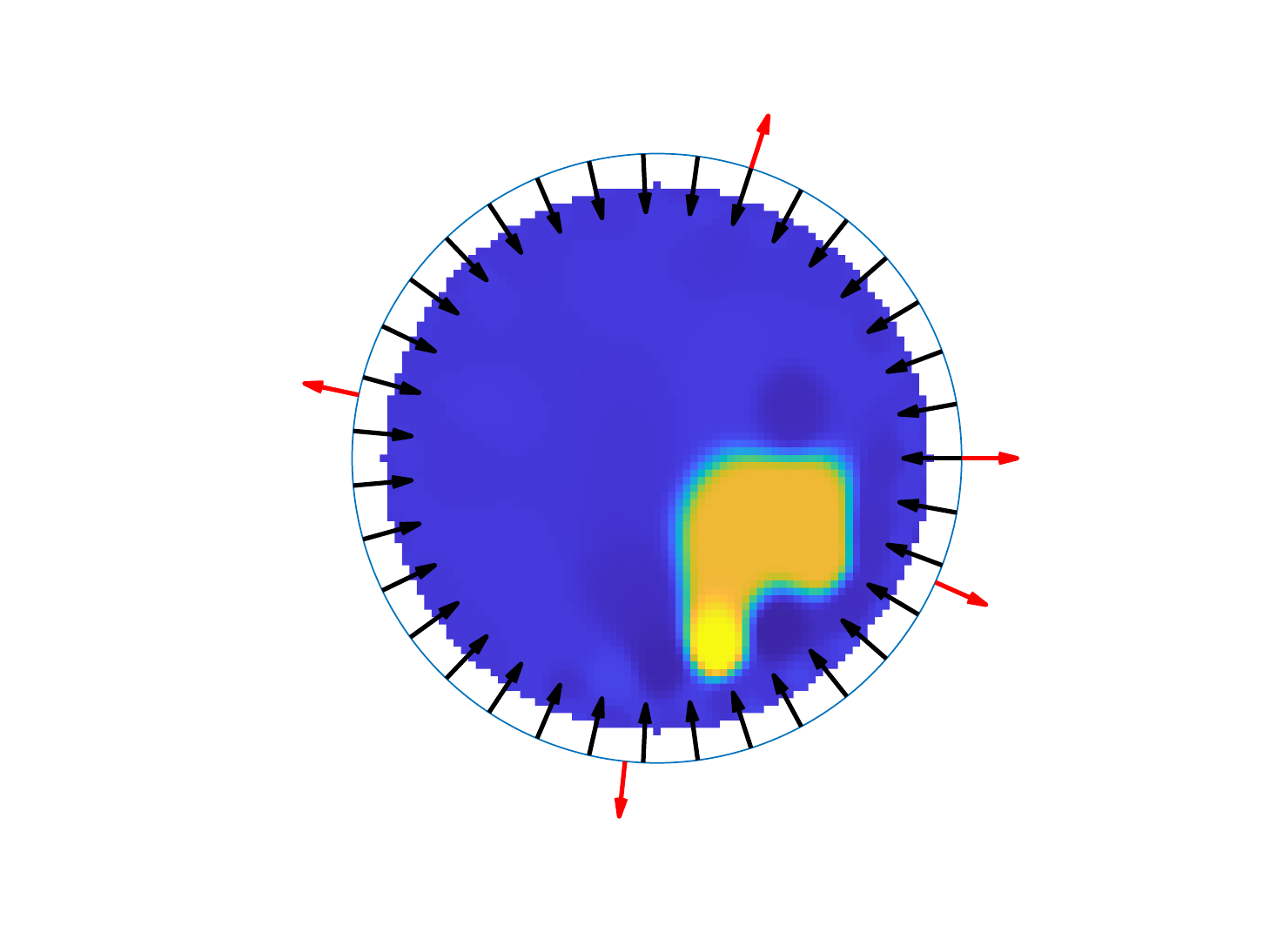}
  \includegraphics[width = 0.3\columnwidth]{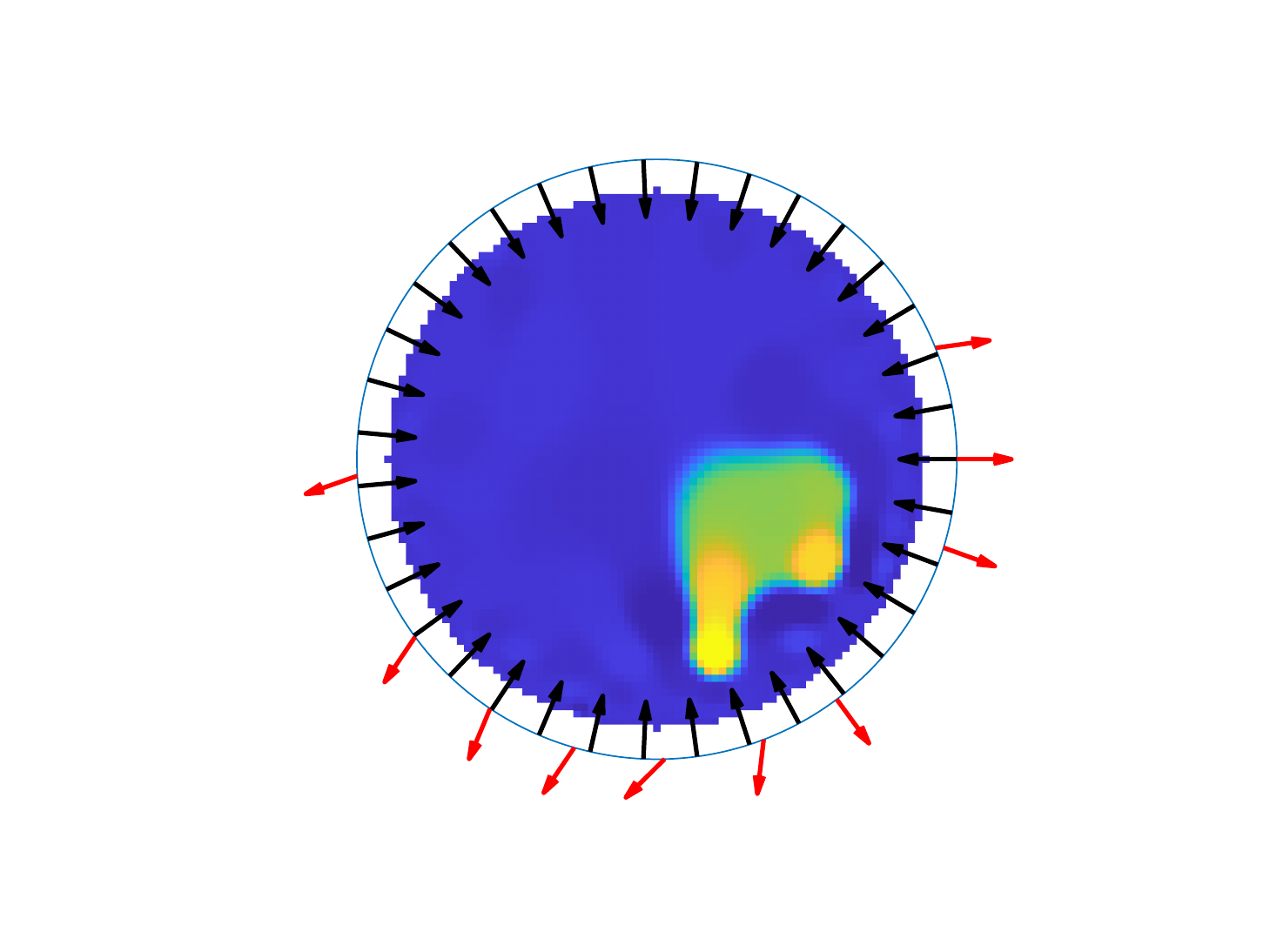}
  \includegraphics[width = 0.3\columnwidth]{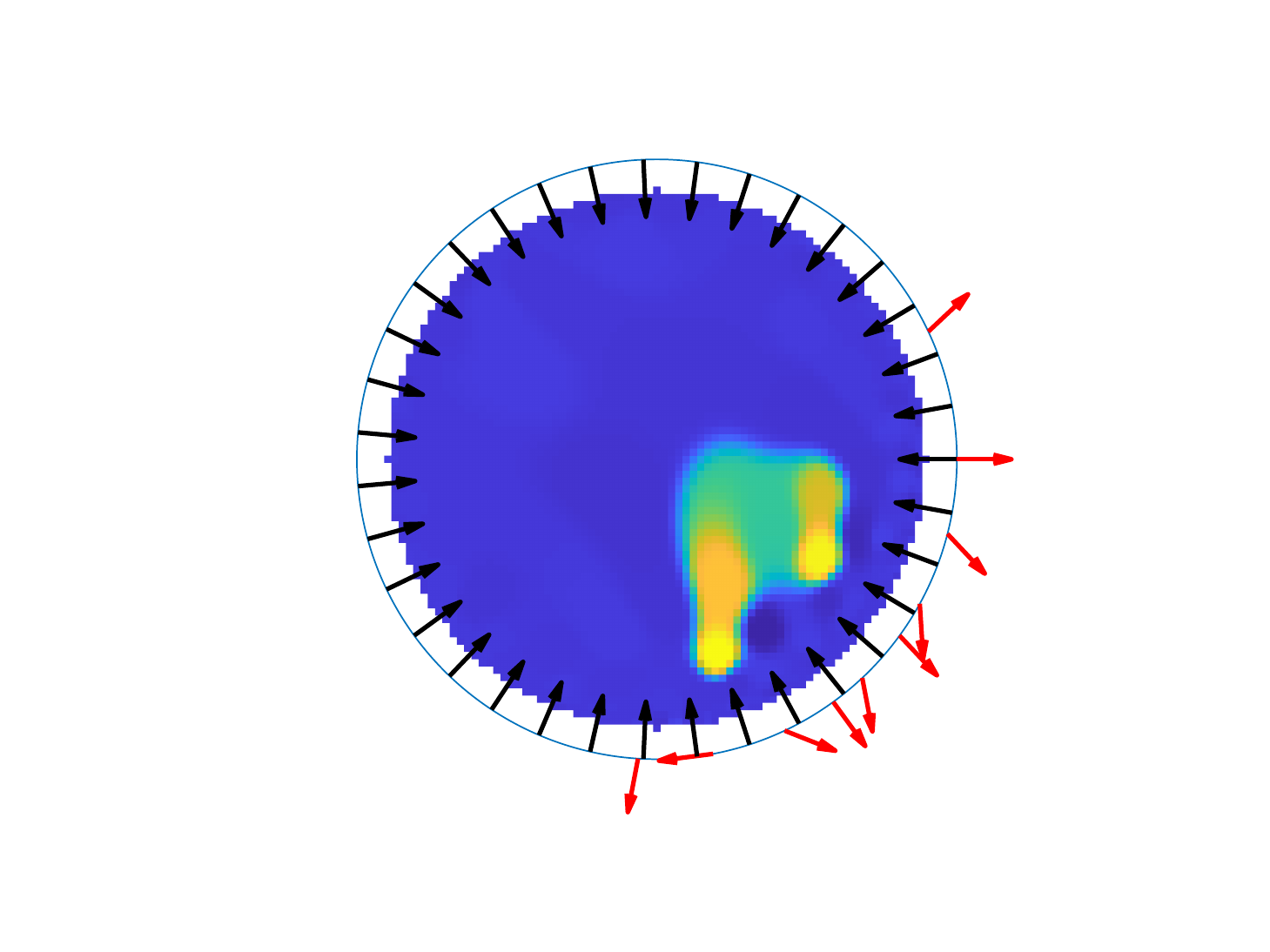}
  \includegraphics[width = 0.3\columnwidth]{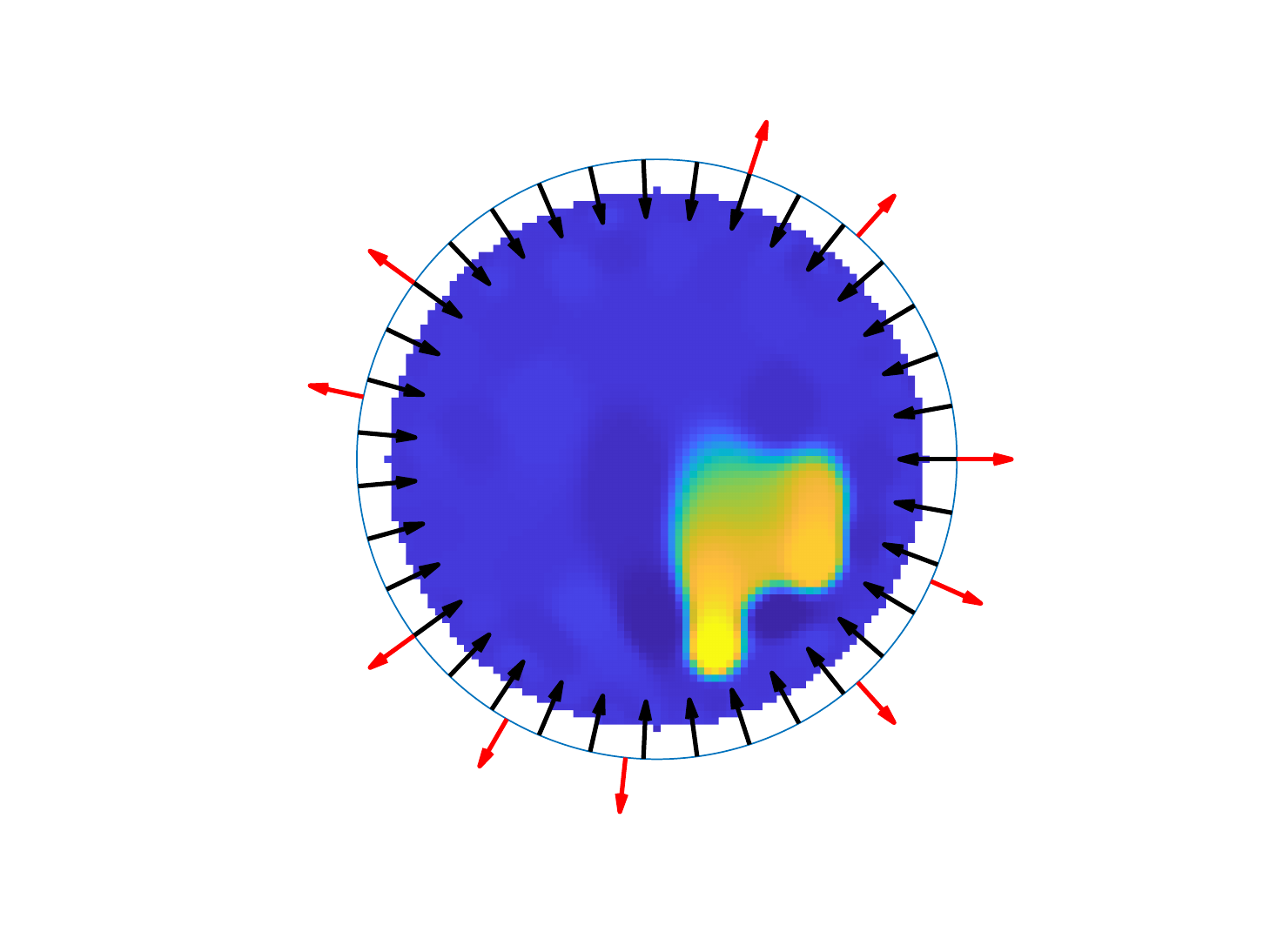}
  \includegraphics[width = 0.3\columnwidth]{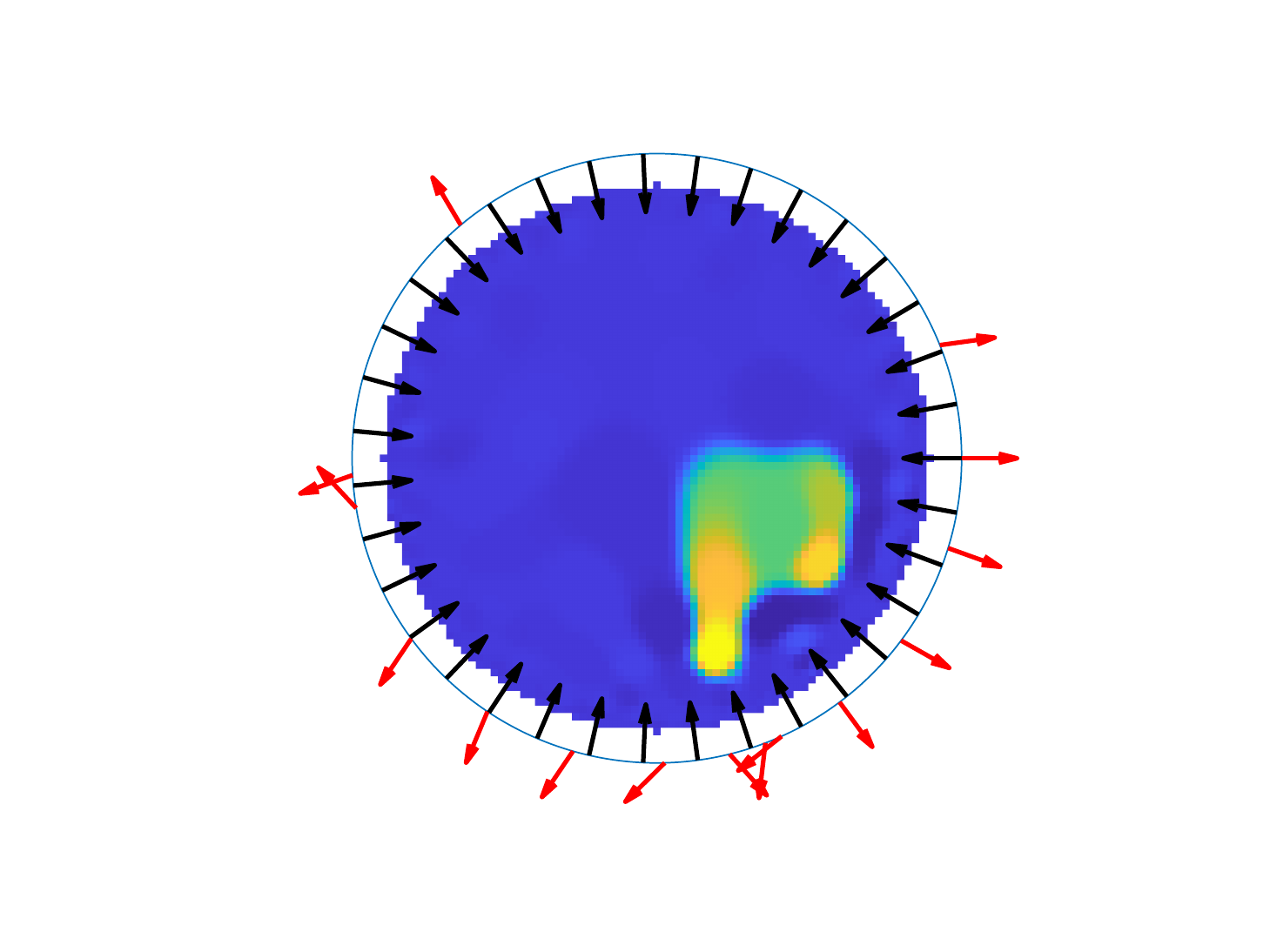}
  \includegraphics[width = 0.3\columnwidth]{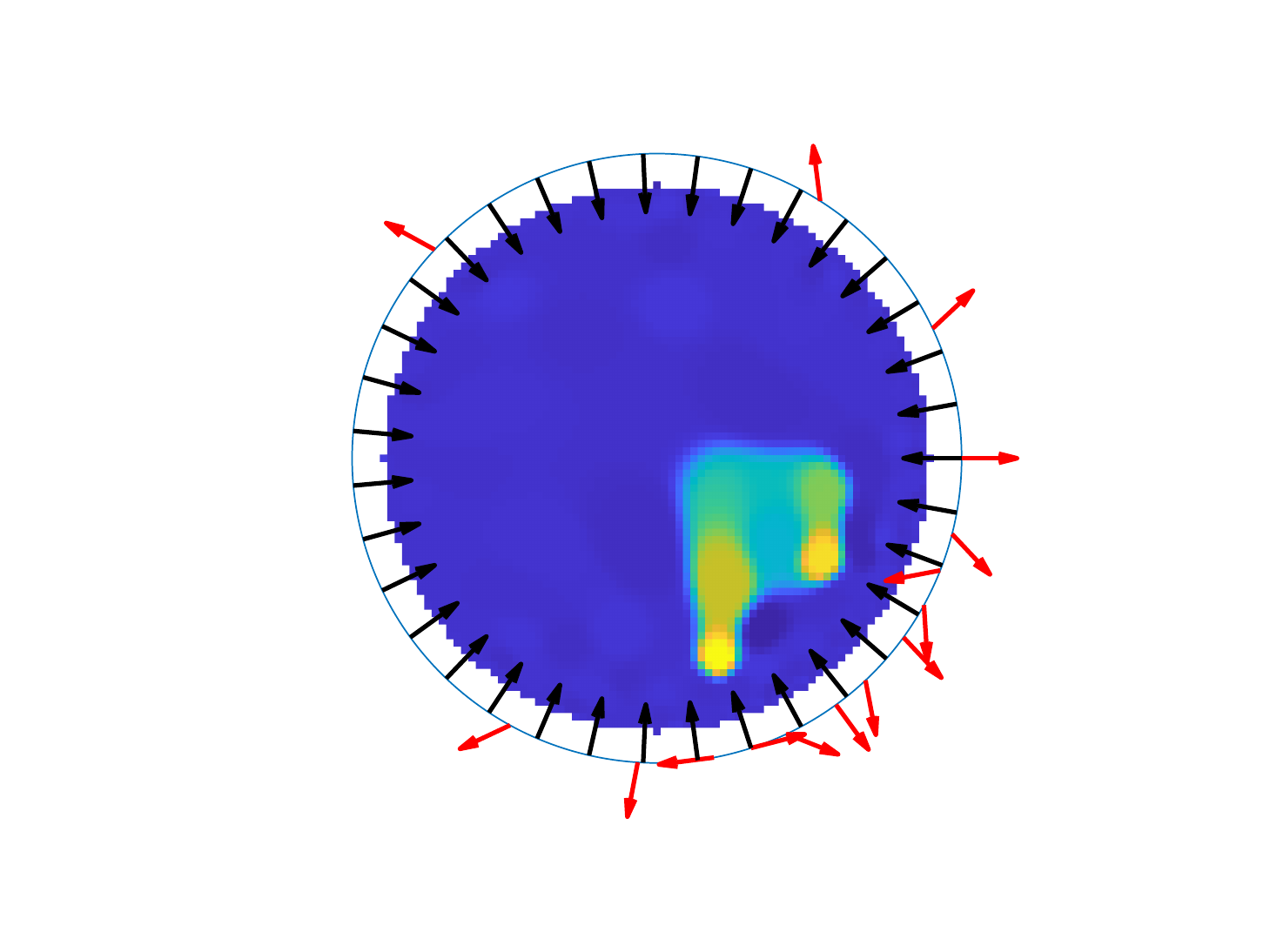}
  \includegraphics[width = 0.3\columnwidth]{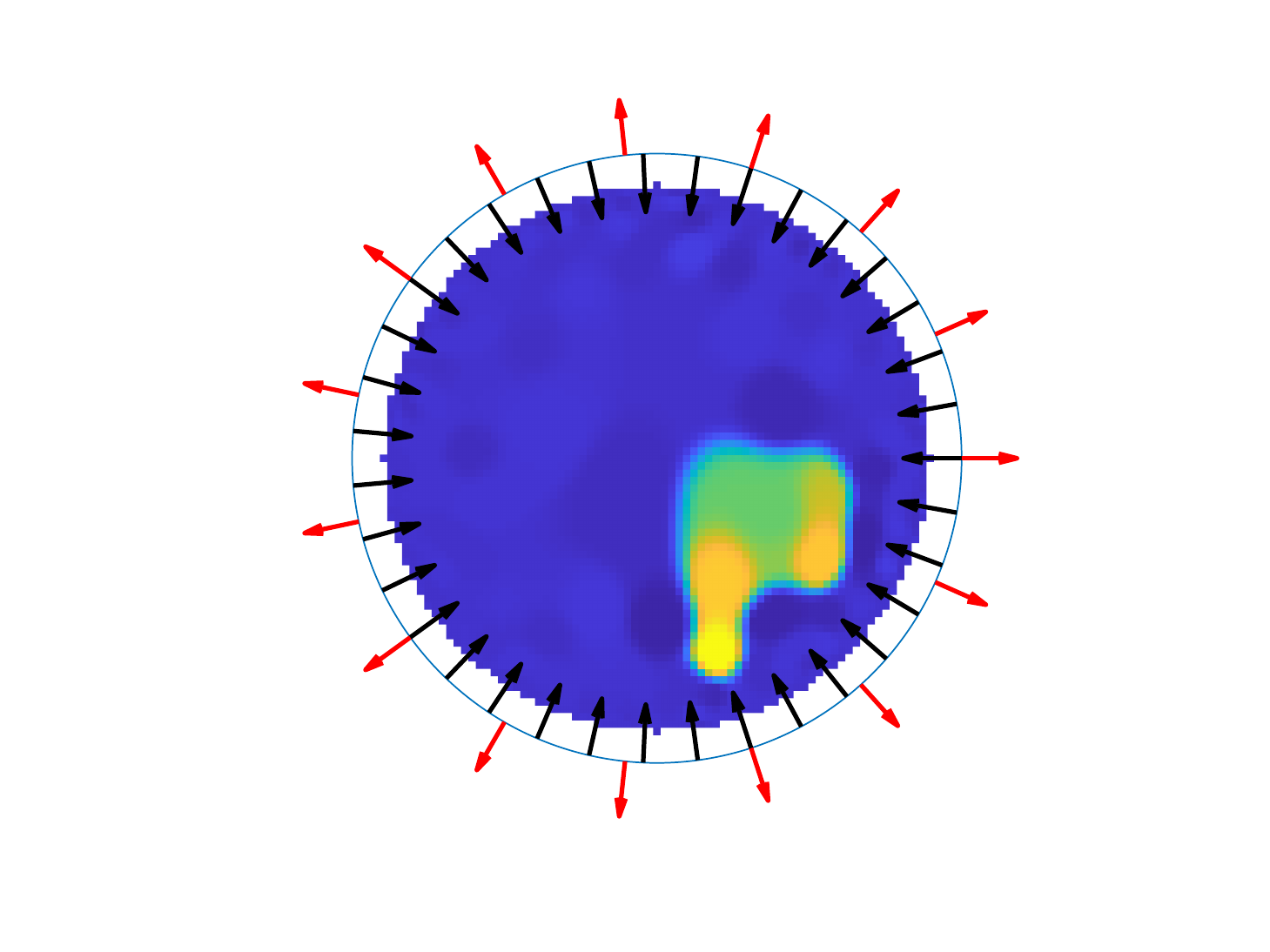}
	\caption{{\sc TV~Test~1.} Reconstructions and activation designs produced by Algorithm~\ref{alg:basic_optimization}. The rows correspond to 5, 10, and 15 activations. Left: Newton's method. Center: exhaustive search. Right: reference designs.}
	\label{fig:TV_test1_2}
\end{figure}

Algorithm~\ref{alg:basic_optimization} is run for a total of $N_a = 15$ iterations, considering all three options for sequentially deducing the specifications of the activations: gradient descent, Newton's method and an exhaustive search. For comparison, Algorithm~\ref{alg:basic_optimization} is also run so that the optimization step on the first line of the for-loop is skipped and the next activation is chosen deterministically from a reference set: The full reference set of activations consists of 15 equiangularly positioned dipoles pointing away from the center of $\Omega$. The activations are introduced one by one into the lagged diffusivity-based  reconstruction algorithm so that new activations are placed in turns onto each quadrant of the measurement circle. The images in the right-hand column of Figure~\ref{fig:TV_test1_2} show the reference setups with five, ten and fifteen activations included. The reference activations are also used as the initial guesses when Newton's method and gradient descent are used for sequential optimization of the activations in  Algorithm~\ref{alg:basic_optimization}.

The right-hand image in Figure~\ref{fig:TV_test1_1} shows the relative $L^2(\Omega)$ reconstruction error as a function of the number of activations for the four options of sequentially choosing the activation designs in Algorithm~\ref{alg:basic_optimization}. It is obvious that the deterministic reference designs perform the worst and the exhaustive search performs in the end the best in terms of the reconstruction error. A bit surprisingly, gradient descent is in this case the second best optimization routine, resulting in even smaller reconstruction errors than the exhaustive search between three and nine activations. This highlights two aspects: (i) The relative performance of gradient descent and Newton's method is highly case-dependent due to the existence of several local minima in the target function; the tendency of these optimization methods to jump between local minima could be tuned by the initial step size parameter $\lambda$ in Algorithms~\ref{alg:GD} and \ref{alg:newton}. (ii) Due to the sequential nature of Algorithm~\ref{alg:basic_optimization}, it is possible that not finding the global, but only local minima by a differentiation-based method when determining individual activations leads to a more optimal {\em combination of activations} than employing the exhaustive search.

According to our experience from other numerical tests not documented here, gradient descent produces on average approximately as optimal designs as Newton's method, but the run time of our nonoptimized MATLAB implementation of Algorithm~\ref{alg:basic_optimization} is typically slightly shorter with Newton's method. This advantage is explained by the possibility to explicitly calculate all needed second derivatives, the small size of the Hessian in Algorithm~\ref{alg:newton} and the better theoretical convergence rate of Newton's method close to a local minimum. In the presented example, Algorithm~\ref{alg:basic_optimization} was 32\% faster with Newton's method than with gradient descent.

\begin{figure}
	\centering
  \includegraphics[width = 0.49\columnwidth]{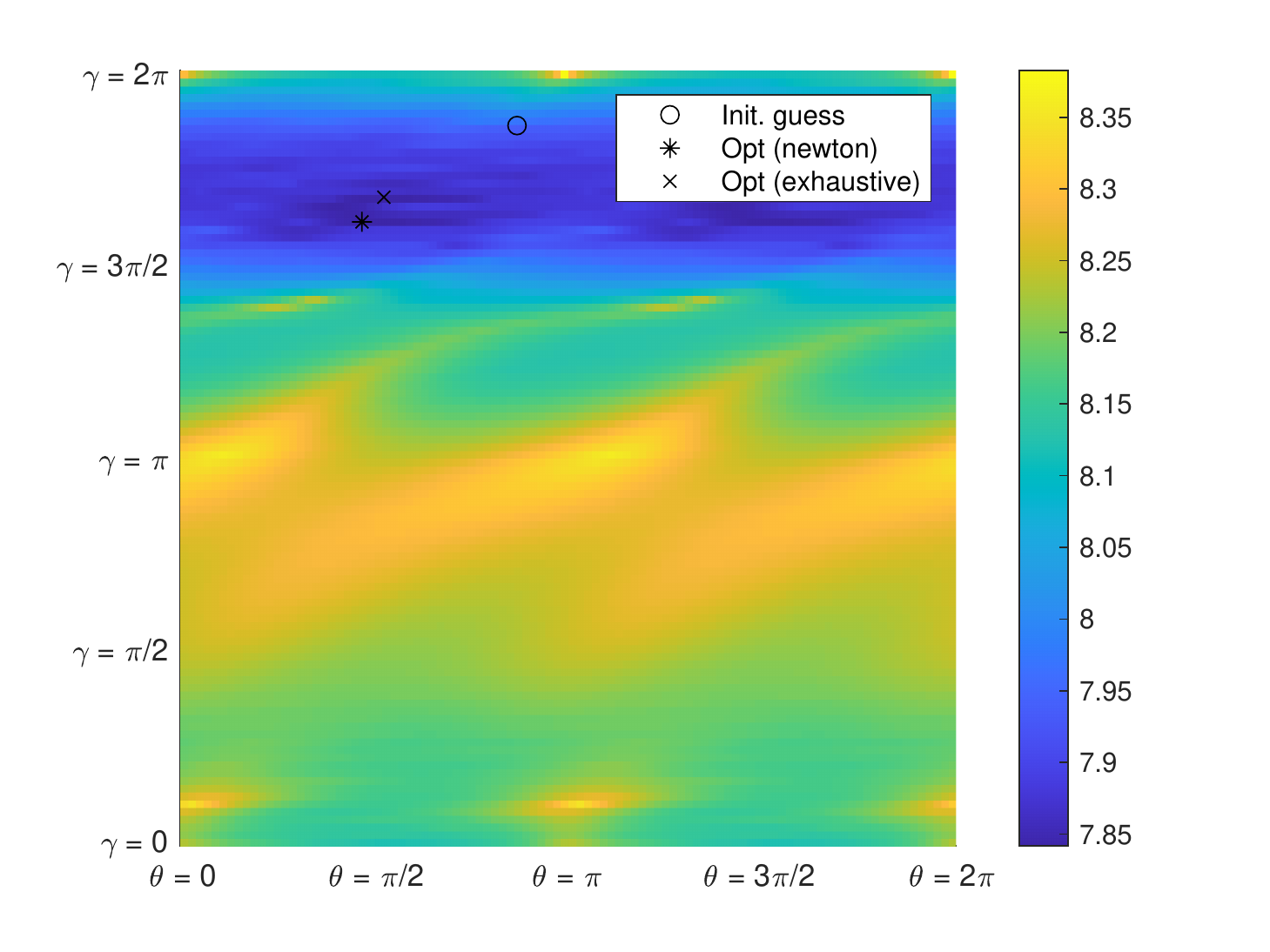}
  \includegraphics[width = 0.49\columnwidth]{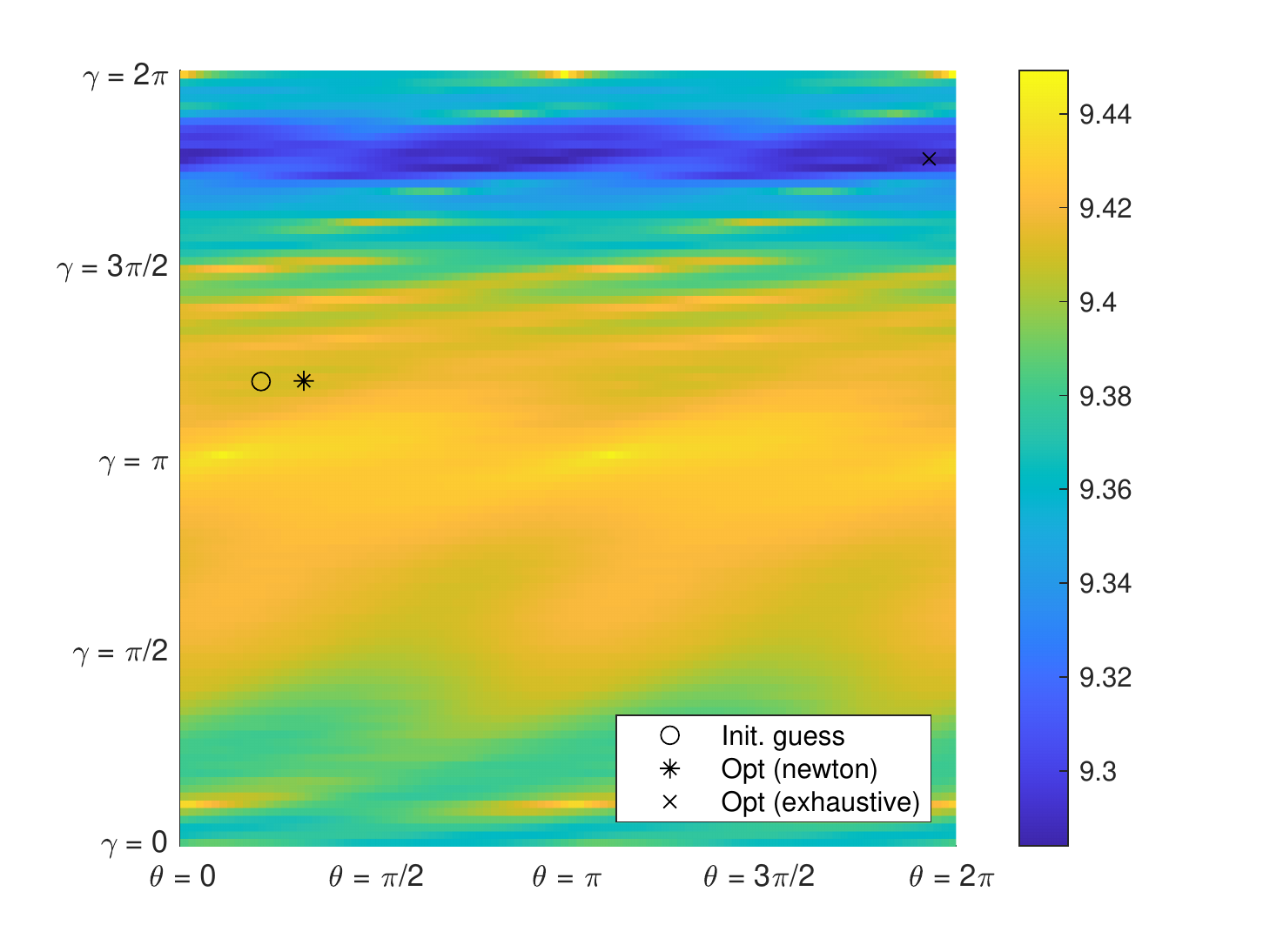}
	\caption{{\sc TV~Test~1.} Visualizations of the A-optimality target $\Psi_{\rm A}$ as a function of the two angular decision variable on $[0, 2\pi] \times [0, 2\pi]$ during a run of Algorithm~\ref{alg:basic_optimization} with Newton's method for the phantom in the left-hand image of Figure~\ref{fig:TV_test1_1}. The vertical axis corresponds to the position of the activation on the measurement circle of radius $1.1\rho$ and the horizontal axis to the direction of the activation dipole. Left: $\Psi_{\rm A}$ after four activations. Right: $\Psi_{\rm A}$ after nine activations. The circle depicts the initial guess for deducing the next activation, the star shows the local minimum found by Newton's method, and the cross indicates the global minimizer.}
	\label{fig:TV_test1_3}
\end{figure}

Figure~\ref{fig:TV_test1_2} shows the reconstructions and optimal designs produced by Algorithm~\ref{alg:basic_optimization} after five, ten and fifteen activations when the sequential designs are deduced by Newton's method (left), exhaustive search (middle) and the deterministic reference procedure (right). According to a visual investigation, the differences in the reconstruction quality are not huge. On the other hand, it is obvious that Newton's method and the exhaustive search find different (sequential) local minima in the optimization target: both methods concentrate the activations close to the location of the P-shaped inclusion, but already after five activations the fine details of the designs are considerably different. This is due to the existence of multiple local minima that makes the performance of differentiation-based optimization routines suboptimal without further modifications, such as using several different initializations. Figure~ \ref{fig:TV_test1_3} illustrates this problem by visualizing the A-optimality target function after four and nine activations when Newton's method is employed. On the left, the optimized fifth activation found by Newton's method is close to the global (sequential) minimizer. However, on the right in case of the tenth activation, we can see that a bad initial guess has resulted in Newton's method getting stuck in a local minimum far from the global optimum.

\subsubsection{TV Test~2: P-shaped inclusion a with larger domain}

\begin{figure}
	\centering
  \includegraphics[scale = 0.5]{figures/tv/target_P.pdf}
  \hspace{-1.5cm}
  \includegraphics[scale = 0.05]{figures/tv/target_P.pdf}
  \hspace{0.5cm}
  \includegraphics[width = 0.56\columnwidth]{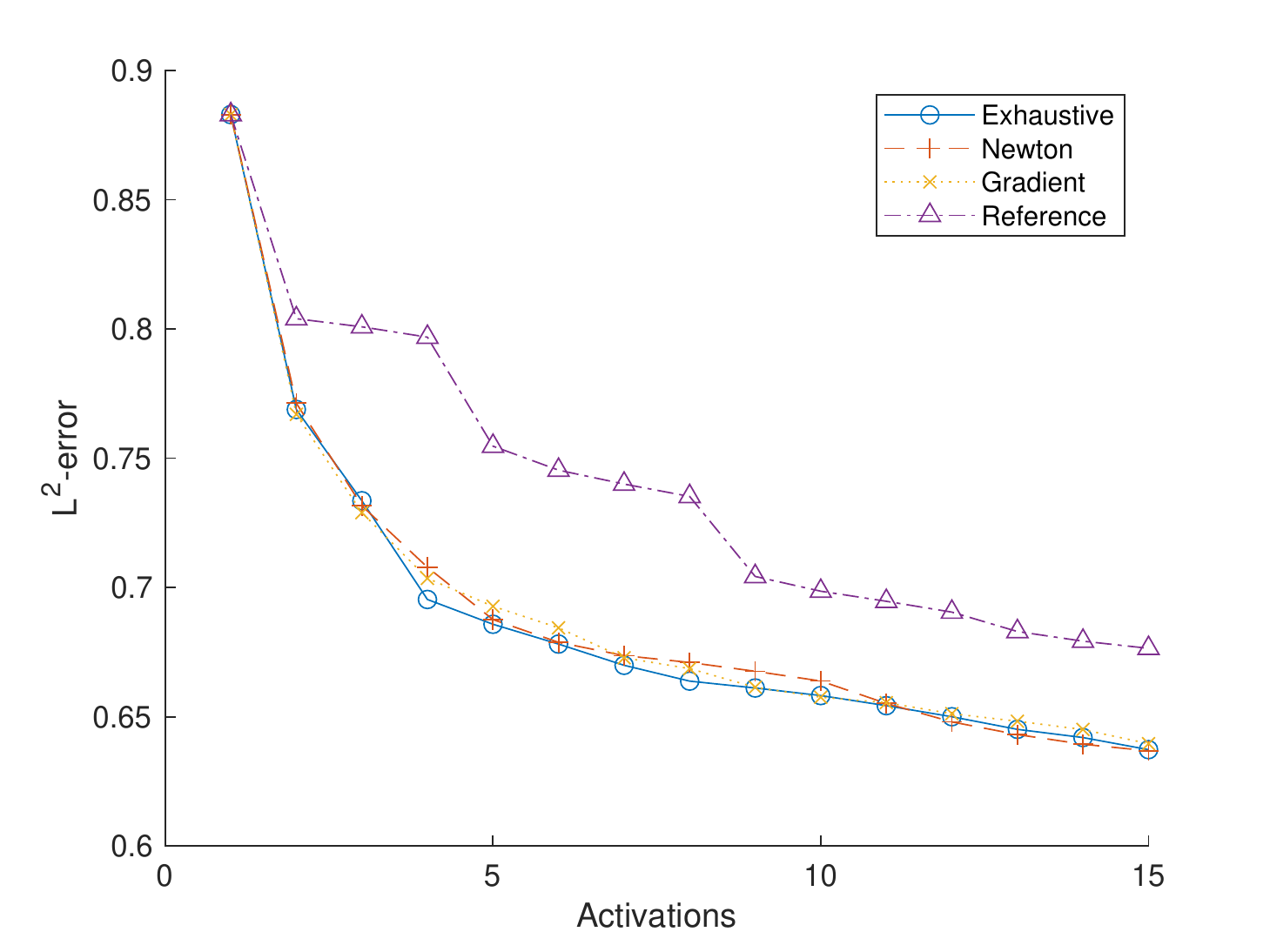}
	\caption{{\sc TV~Test~2.} Left: target MNP concentration compared to the much smaller target from Figure~\ref{fig:TV_test1_1}. Right: relative $L^2(\Omega)$ reconstruction errors as functions of the number of activations in Algorithm~\ref{alg:basic_optimization}  when the specifications of the activations are deduced by Newton's method, Gradient descent, exhaustive search and the deterministic reference procedure.}
	\label{fig:TV_test2_1}
\end{figure}

Our second test on sequential edge-promoting Bayesian ODE essentially repeats the previous experiment with the noise and size parameters reset to $\eta = 0.01$ and $\rho = 5$,~i.e.,~the noise level is lower  but the object is larger than previously. The right-hand image of Figure~\ref{fig:TV_test2_1} shows the relative $L^2(\Omega)$ reconstructions errors for the four methods for determining the specifications of the activation dipoles in Algorithm~\ref{alg:basic_optimization}. 
This time, Newton's method and gradient descent lead to roughly equally accurate reconstructions as the exhaustive search, with the deterministic procedure for choosing the activation designs again performing the worst. In fact, compared to the previous test with a smaller $\Omega$, the deterministic reference procedure results now in much larger reconstruction errors compared to the other three options, presumably due to the larger size of the imaged object that makes activations,~e.g.,~on the top-left side of $\Omega$ less informative. 

\begin{figure}
	\centering
  \includegraphics[width = 0.3\columnwidth]{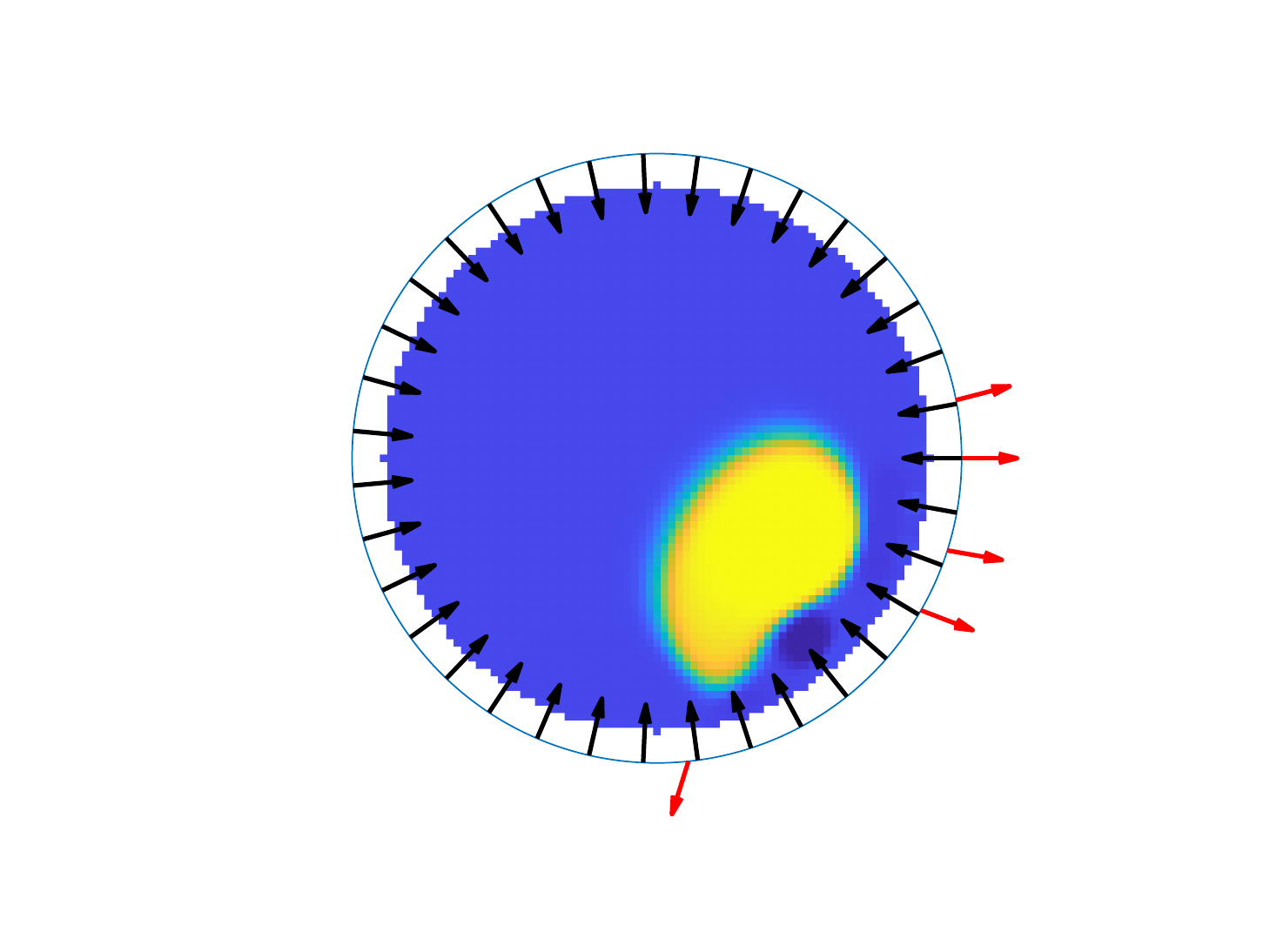}
  \includegraphics[width = 0.3\columnwidth]{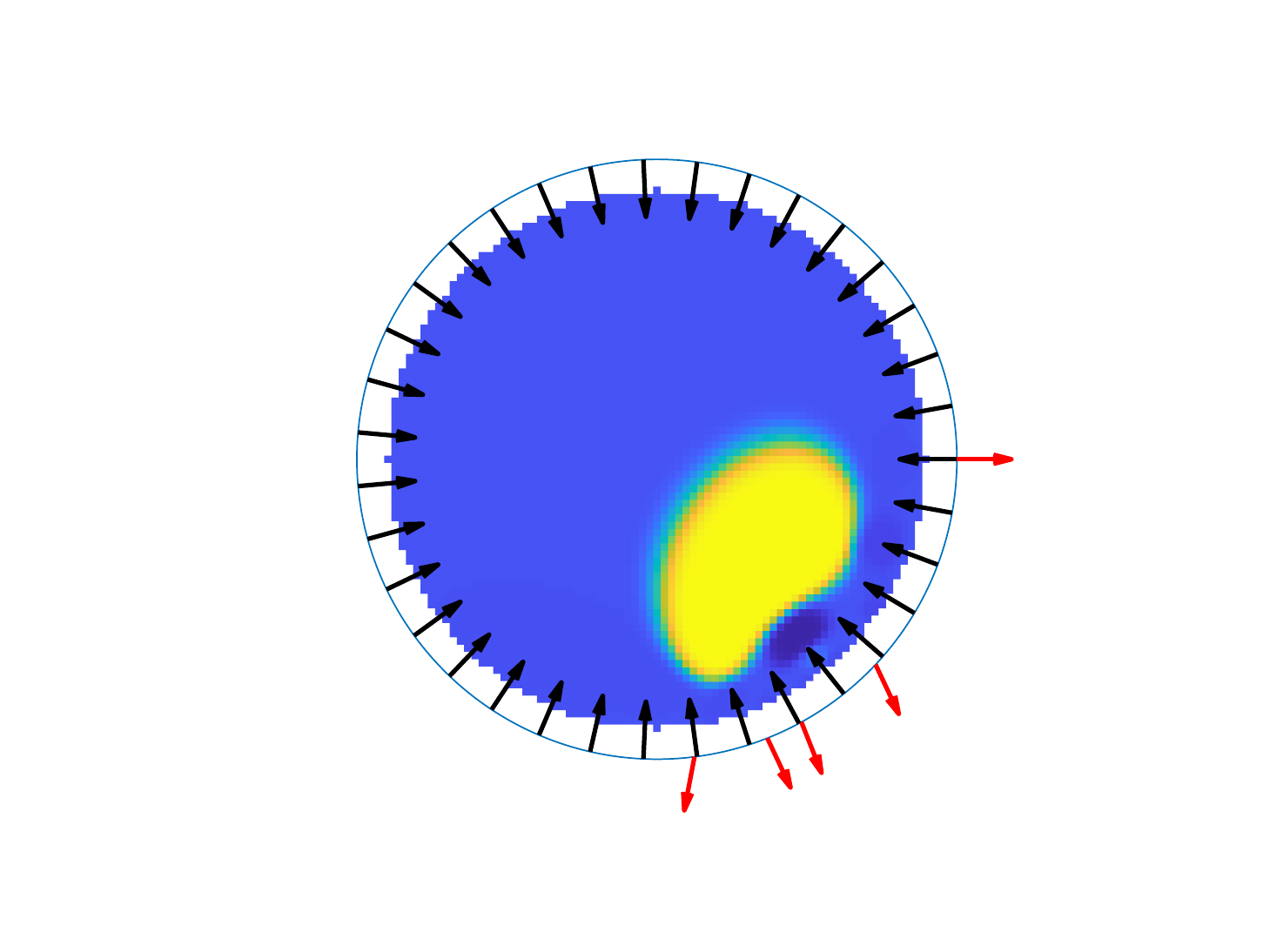}
  \includegraphics[width = 0.3\columnwidth]{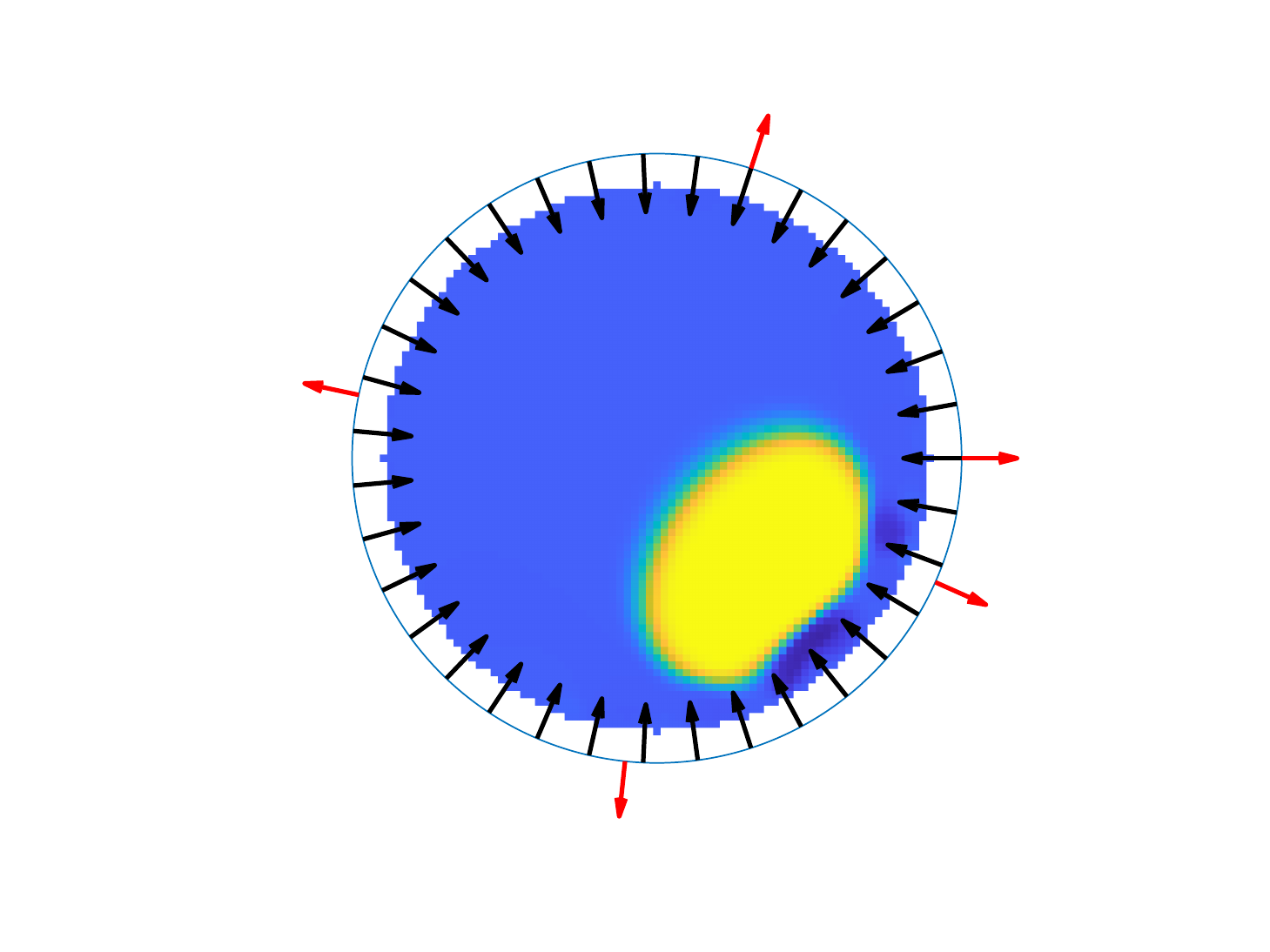}
  \includegraphics[width = 0.3\columnwidth]{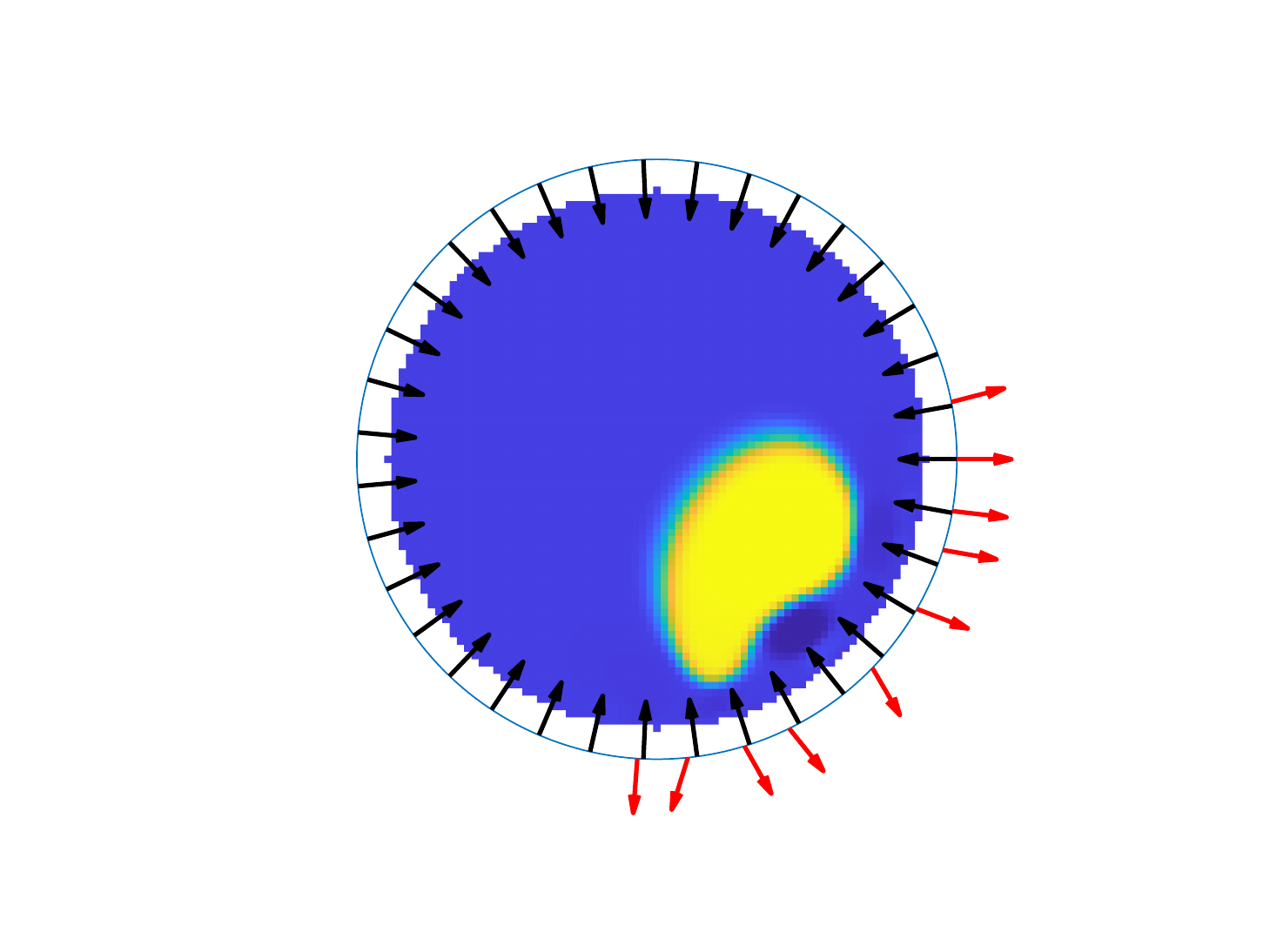}
  \includegraphics[width = 0.3\columnwidth]{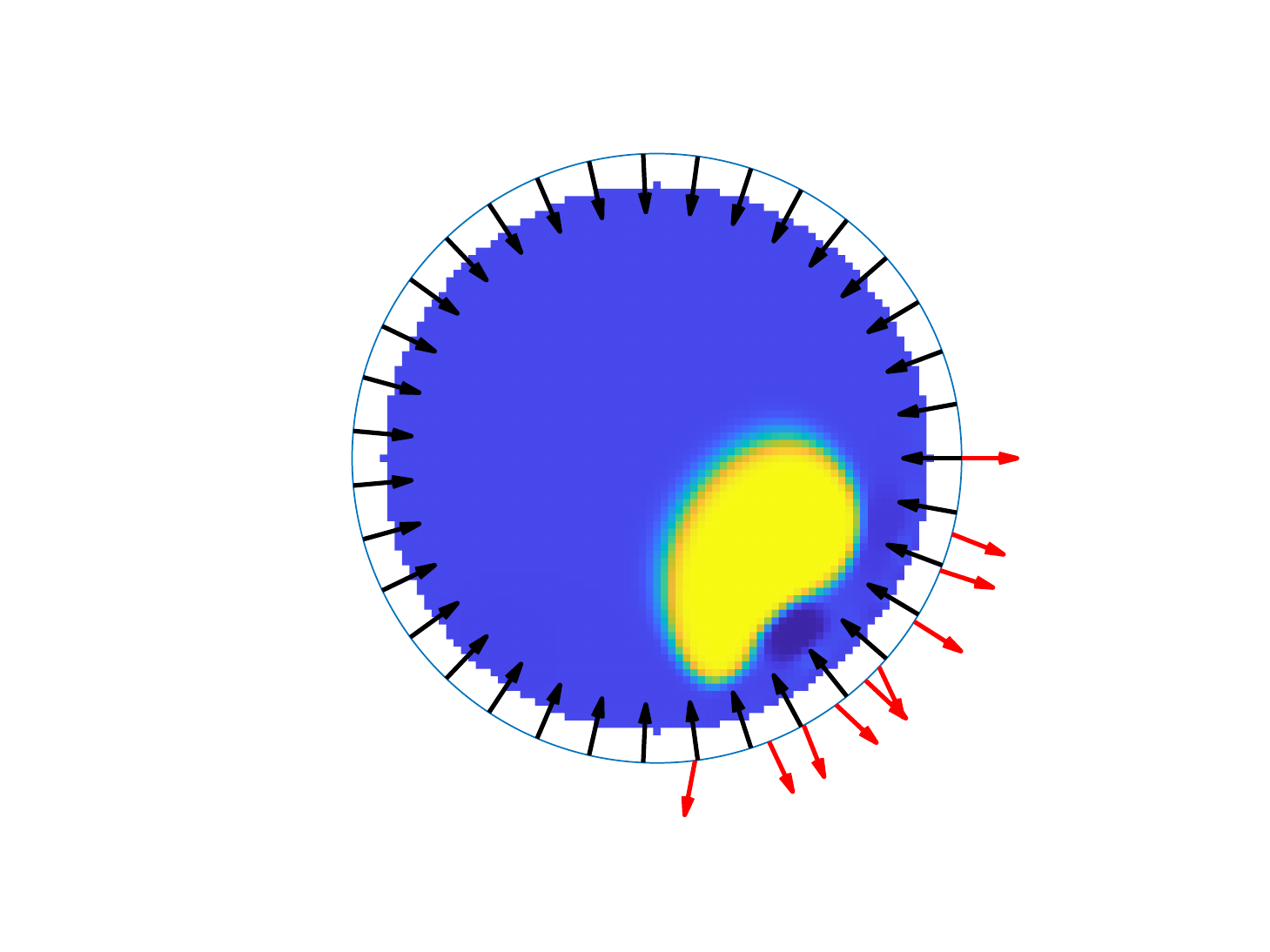}
  \includegraphics[width = 0.3\columnwidth]{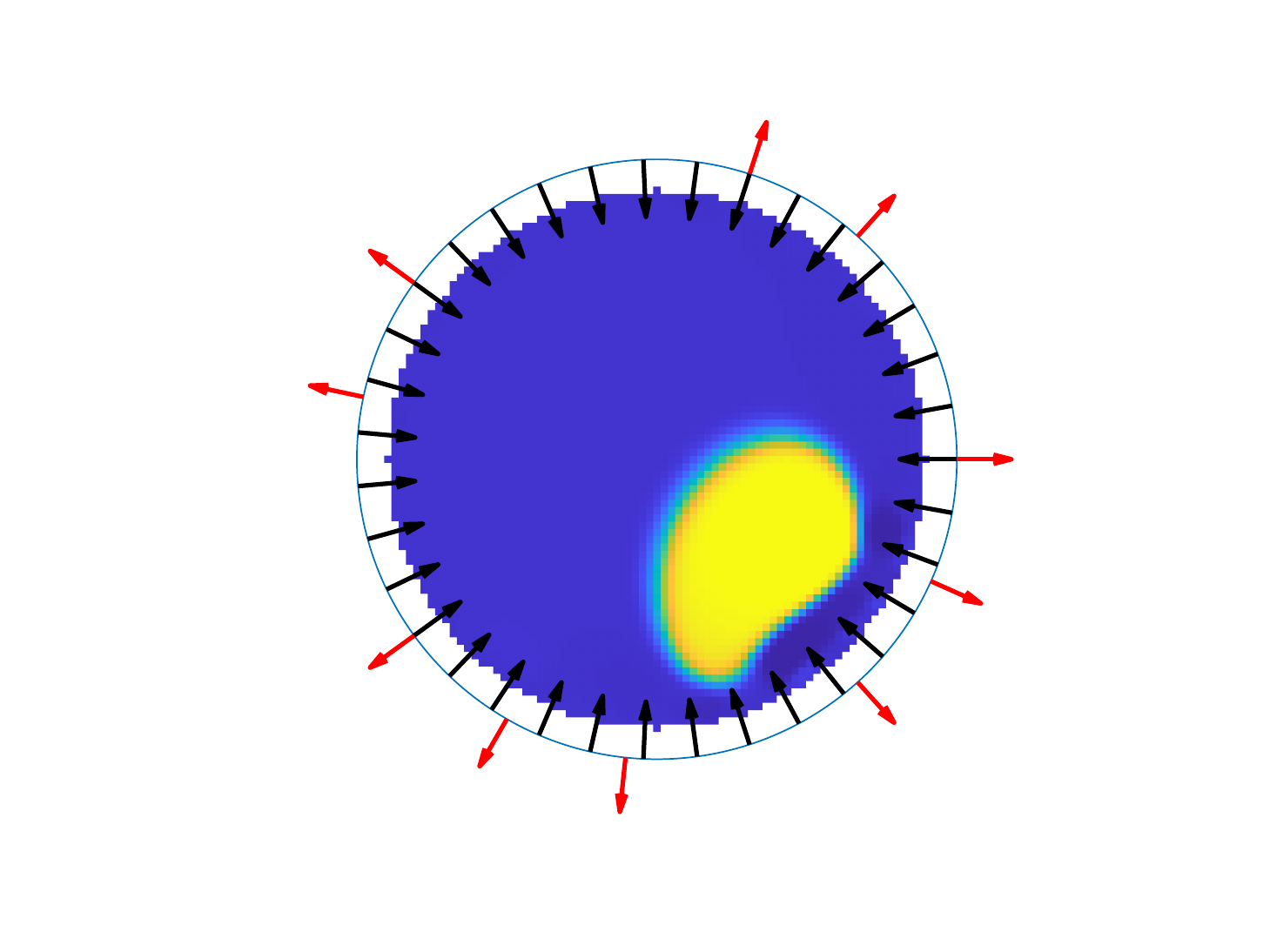}
  \includegraphics[width = 0.3\columnwidth]{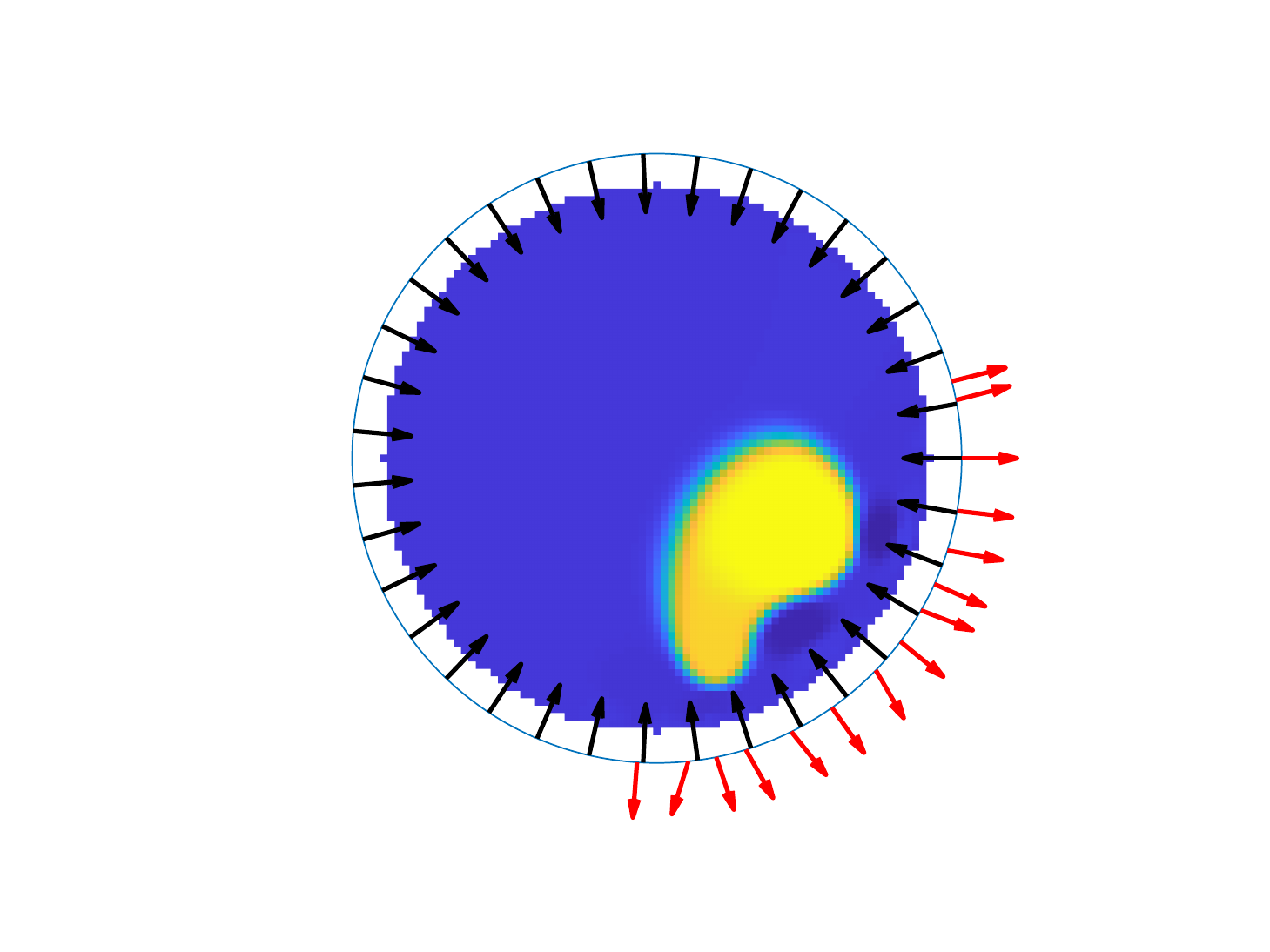}
  \includegraphics[width = 0.3\columnwidth]{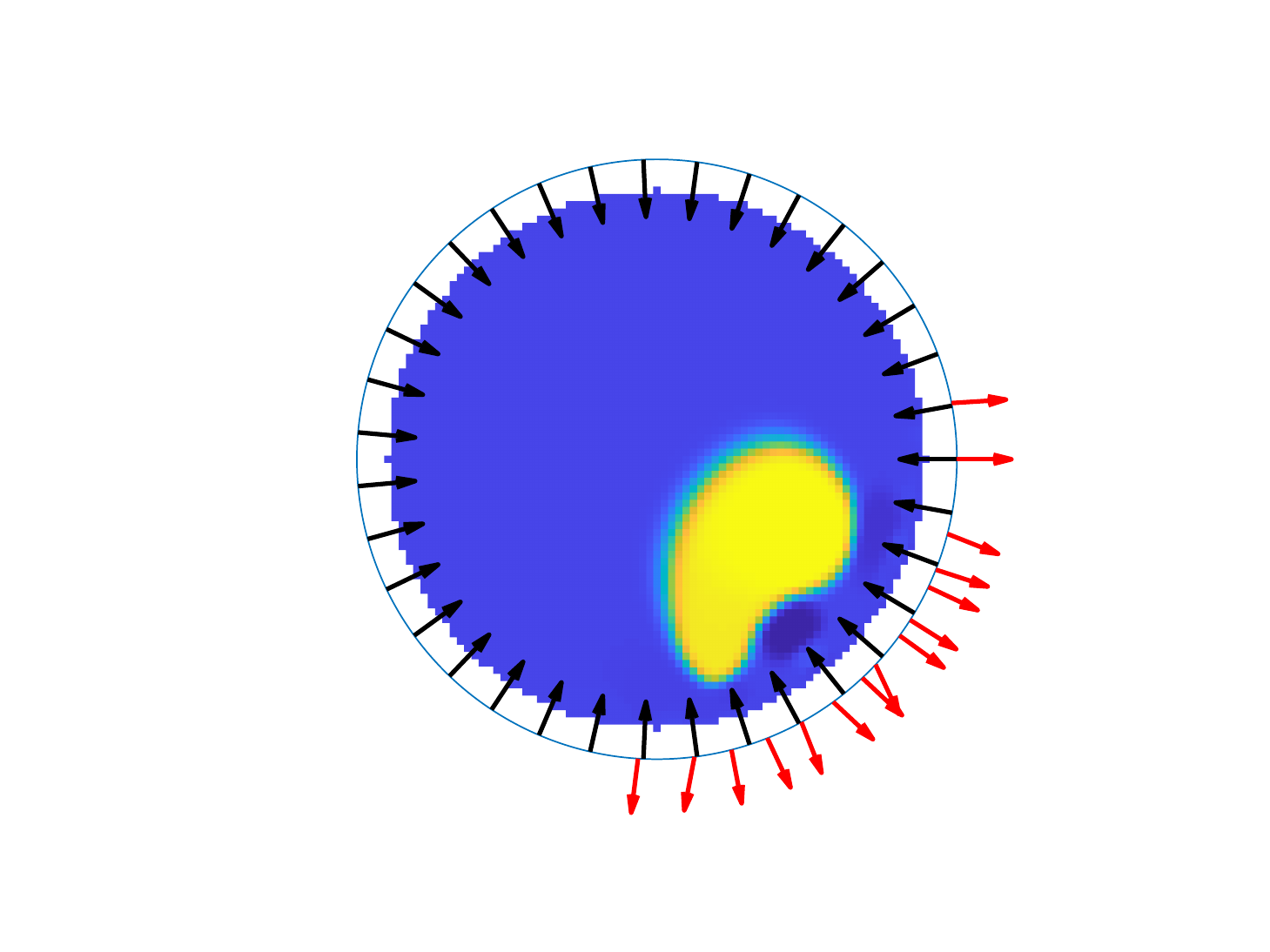}
  \includegraphics[width = 0.3\columnwidth]{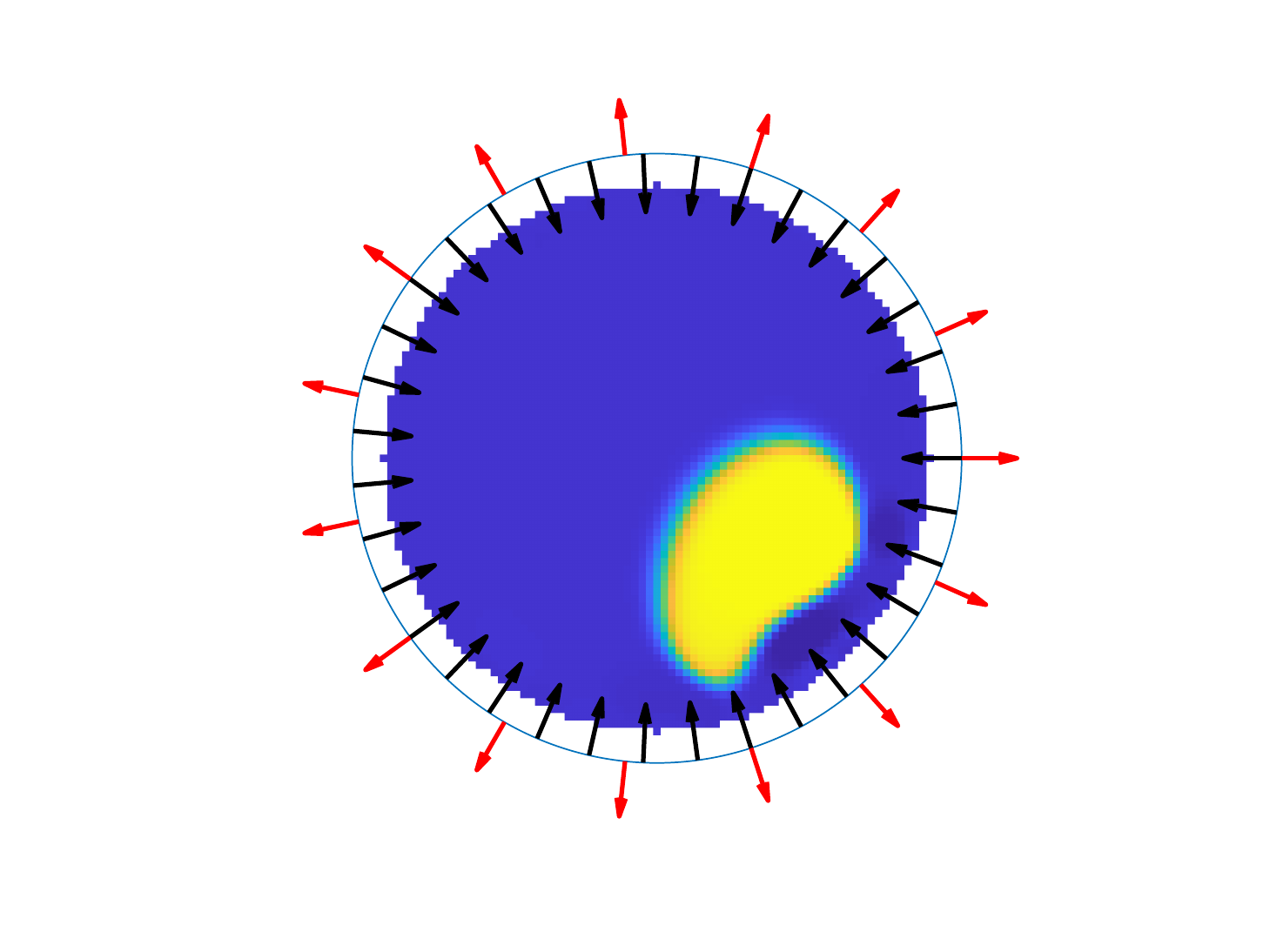}
	\caption{{\sc TV~Test~2.} Reconstructions and activation designs produced by Algorithm~\ref{alg:basic_optimization}. The rows correspond to 5, 10, and 15 activations. Left: Newton's method. Center: exhaustive search. Right: reference designs.}
	\label{fig:TV_test2_2}
\end{figure}

Figure~\ref{fig:TV_test2_2}, which is organized in the same way as Figure~\ref{fig:TV_test1_2} for the previous test, demonstrates that the optimal designs produced by Algorithm~\ref{alg:basic_optimization} with both Newton's method and the exhaustive search have the activations more tightly packed around the interesting area in $\Omega$ than in the preceding test. Although these optimized designs seem qualitatively very similar, there are differences in their fine details due to the effect of local minima. As the difference between the reference and optimized designs is considerable, it is hardly surprising that there is also a larger difference between the corresponding reconstruction errors compared to the previous test with a smaller target domain.

\section{Concluding remarks}
\label{sec:conclusion}
This work investigated the application of Bayesian OED techniques to MRXI in a simplified two-dimensional simulated setting. The numerical experiments tested the general applicability of the studied approach to simultaneous optimization of several activations with a Gaussian prior and to sequentially choosing optimal activations with a smoothened TV prior following the ideas in \cite{Helin22}. The presented results demonstrate that choosing the activations according to the A-optimality criterion can improve reconstruction quality, especially when the imaged object or the ROI is nonsymmetric, which makes symmetric reference designs suboptimal.

The numerous local minima in the optimization targets (cf.~Figure~\ref{fig:TV_test1_3}) constitute a major obstacle for efficiently applying Bayesian OED to realistic MRXI, a problem for which no solutions were proposed in this work.  
Moreover, no effort was made to model practical measurement settings of MRXI, that is, to account for the inherently three-dimensional nature of the measurements, the shape of the imaged object, accurate models for activations and measurements, and the effect of realistic values for the involved physical quantities. Designing optimization algorithms that are robust to local minima and working with more detailed and realistic measurement models are interesting topics for future research on Bayesian OED for MRXI.

\section*{Acknowledgments}
We would like to thank Tommi Huhtinen and Laura Rautiainen for performing initial computational tests on Bayesian OED for MRXI during their summer internships at the Department of Mathematics and Systems Analysis, Aalto University in 2022.

\bibliographystyle{acm}
\bibliography{mrxi_opt_refs}
\end{document}